\documentclass[11pt]{amsart}
\usepackage{amsmath,amssymb,amsfonts,url,mathptmx}
\usepackage{color}
\usepackage{appendix}
\numberwithin{equation}{section}

\newtheorem{thm}{Theorem}[section]
\newtheorem{lem}{Lemma}[section]
\newtheorem{rem}{Remark}[section]

\newcommand{\hdot}{^\text{\r{}}\hspace{-.33cm}H}
\usepackage[marginal]{footmisc}

\begin{document}
\title[$SU(3)$ Toda system]{Estimates of bubbling solutions of $SU(3)$ Toda systems at critical parameters-Part 1}
\keywords{}

\author{Lina Wu}\footnote{Lina Wu is partially supported by the China Scholarship Council (No.201806210165).}

\author{Lei Zhang}\footnote{Lei Zhang is partially supported by a Simons Foundation Collaboration Grant}

\address{Department of Mathematical Sciences \\
         Tsinghua University \\
         No.1 Qinghuayuan, Haidian District,
         Beijing China, 100084 }
\email{wln16@mails.tsinghua.edu.cn}

\address{Department of Mathematics\\
        University of Florida\\
        1400 Stadium Rd\\
        Gainesville FL 32611}
\email{leizhang@ufl.edu}

\date{\today}

\begin{abstract}
For regular $SU(3)$ Toda systems defined on Riemann surface, we initiate the study of bubbling solutions if parameters $(\rho_1^k,\rho_2^k)$ are both tending to critical positions: $(\rho_1^k,\rho_2^k)\to (4\pi, 4\pi N)$ or $(4\pi N, 4\pi)$ ($N>0$ is an integer). We prove that there are at most three formations of bubbling profiles, and for each formation  we identify leading terms of $\rho_1^k-4\pi$ and $\rho_2^k-4\pi N$, locations of blowup points and comparison of bubbling heights with sharp precision. The results of this article will be used as substantial tools for a number of degree counting theorems, critical point at infinity theory in the future.
\end{abstract}


\maketitle

\section{Introduction} It is well known that a lot of models in various disciplines of sciences are described by second order elliptic systems with exponential nonlinear terms.  Among all such systems and equations defined in two dimensional spaces,  Toda system is probably the most important one, as it has profound connections with nonabelian gauge field in Chern-Simon's gauge theory
(see \cite{dunne,han-yang,nie,yang1} etc), the solution of its reduced form represents metric with prescribed Gauss curvature with conic singularities (see \cite{bat-mal,chani-kiess,fang-lai,jost-wang,troy}, etc). In this article we start from its simplest form: $SU(3)$Toda system and study the behavior of bubbling solutions in a critical situation.

Let $(M,g)$ be a Riemann surface with $vol(M)=1$ for convenience. The main equation we consider in this article is
\begin{equation}\label{main-eq-1}
\left\{\begin{array}{ll}
\Delta_g \tilde u_1+2\rho_1(\frac{h_1e^{\tilde u_1}}{\int_M h_1e^{\tilde u_1}}-1)-\rho_2(\frac{h_2e^{\tilde u_2}}{\int_M h_2e^{\tilde u_2}}-1)=0,\\
\\
\Delta_g \tilde u_2-\rho_1(\frac{h_1e^{\tilde u_1}}{\int_M h_1e^{\tilde u_1}}-1)+2\rho_2(\frac{h_2e^{\tilde u_2}}{\int_M h_2e^{\tilde u_2}}-1)=0.
\end{array}
\right.
\end{equation}
where $\rho_1,\rho_2>0$ are constants, $h_1$, $h_2$ are positive, smooth functions on $M$, $\Delta_g$ is Laplace-Beltrami operator $(-\Delta_g\ge 0)$. The coefficient matrix in (\ref{main-eq-1}):
$$\left(\begin{array}{cc}
2 & -1 \\
-1 & 2
\end{array}
\right) $$
is the simplest Cartan matrix, which has its more general form $K$ as
\begin{equation}\label{cartan-m}
\mathbf{K}=\left(\begin{array}{ccccc}
2 & -1 & 0 & ... & 0 \\
-1 & 2 & -1 & ... & 0 \\
\vdots & \vdots &  &  & \vdots\\
0 & ... & -1 & 2 & -1 \\
0 & ... &  0 & -1 & 2
\end{array}
\right ).
\end{equation}
The research of Toda systems has sustained over half an century of intensive and extensive investigations and the new results coming from different groups, new inspirations ignited by various perspectively, still make the whole field as dynamic as ever. It is simply impossible to list all the important references in any reasonable manner, so we focus on its connections with holomorphic curves in $\mathbb C\mathbb P^n$, flat $SU(n+1)$ connection, complete integrability and harmonic sequences etc, the interested readers may look into \cite{bolton-1},\cite{bolton-2},\cite{calabi},\cite{chern},\cite{doliwa},\cite{guest},\cite{leznov},\cite{lwy}, etc for more exquisite presentations in many directions.

The variational form of the more general system of $n-$equations is
$$J_{\rho}(\tilde u)=\frac 12\int_M\sum_{i,j=1}^nk^{ij}\nabla_g\tilde u_i\nabla_g \tilde u_jdV_g-{\sum_{i=1}^n\rho_i\log \int_M h_ie^{\tilde u_i}dV_g,}$$
where $(k^{ij})_{n\times n}$ is the inverse of the Cartan matrix, $\tilde u=(\tilde u^1,...,\tilde u^n)\in \hdot^{1}(M)$, which is defined as
$$ \hdot^{1}(M):=\{v\in L^2(M);\quad \nabla v\in L^2(M), \mbox{and }\,\, \int_M v dV_g=0\}. $$

Note that the integral equal to $0$ assumption is due to that fact that adding any constant vector to a solution of (\ref{main-eq-1}) gives rise to another solution. That said, we normalize the solution and write the system in an equivalent form:
Let
$$u_i=\tilde u_i-\log \int_M h_i e^{\tilde u_i}dV_g,\quad i=1,2, $$
then we have
\begin{equation}\label{normal}
\int_M h_ie^{u_i}=1,\quad i=1,2.
\end{equation}
and the main system becomes
\begin{equation}\label{main-2}
\left\{\begin{array}{ll}
\Delta_g u_1+2\rho_1 h_1 e^{u_1}-\rho_2 h_2 e^{u_2}=2\rho_1-\rho_2,\\
\\
\Delta_g u_2-\rho_1 h_1 e^{u_1}+2\rho_2 h_2 e^{u_2}=2\rho_2-\rho_1.
\end{array}
\right.
\end{equation}
Here we recall that
\begin{equation}\label{h1h2}
h_1,h_2 \mbox{ are positive smooth functions on } M.
\end{equation}
Many times we shall use the Green's function corresponding to $-\Delta_g$:
$$-\Delta_{g} G(p,y)=\delta_p-1, \quad \mbox{in}\quad M, $$
and in a neighborhood of $x$ we shall use
\begin{equation}\label{green-ex}
G(x,y)=-\frac{1}{2\pi}\log |x-y|+\beta(x,y).
\end{equation}

In a recent important work \cite{lwyz-apde} a priori estimates for all singular rank 2 Toda systems have been obtained, which laid down the foundation for a program of degree counting for all rank 2 systems. In particular for $SU(3)$ regular Toda systems, a prior estimate holds as long as $\rho_1$ and $\rho_2$ are not multiples of $4\pi$. This means for $4\pi m<\rho_1<4\pi (m+1)$, $4\pi N<\rho_2<4\pi (N+1)$, there should be a topological degree that only depends on $m,N$ and the genus of the manifold $M$.  On the other hand, it is crucial to understand the asymptotic behavior of bubbling solutions when $(\rho_1,\rho_2)$ tends to the grid points $(4\pi m, 4\pi n)$. A lot of work on blowup analysis and degree counting has been done when \emph{one} of $\rho_i^k$s crosses a multiple of $4\pi$ while the other stays away from $4\pi N$ (see \cite{lee-1,lee-2} for example). In this article we initiate a direct attack to the study of bubbling solutions when both parameters are tending to critical positions. The study of this case is also closely related to theory of ``critical point at infinity" of Bahri-Coron \cite{bahri-1,bahri-2}, the study of which will be carried out in forthcoming works.

Let $u^k=(u_1^k,u_2^k)$ be a sequence of blowup solutions, we use $\lambda_i^k$ to denote the maximum of $u_i^k$ on $M$ for $i=1,2$ and we assume
\begin{equation}\label{both-b}
\lambda_i^k\to \infty, \,\, \mbox{ for }\,\, i=1,2,\quad e^{-\lambda_i^k/4}\lambda_j^k\to 0, \quad i\neq j.
\end{equation}
Note that $u_1^k$ and $u_2^k$ may not tend to infinity at the same point. We say $p$ a blow up point of $u^k$ if there is a sequence of points $x_k\to p$ such that
$$\max\{u_1^k(x_k),u_2^k(x_k)\}\to \infty. $$
The last inequality in (\ref{both-b}) is for convenience and is not a problem in applications. Also throughout the paper we don't distinguish sequences and subsequences.

One main feature of Toda system, as well as many second order elliptic equations/systems with exponential terms, is the concentration phenomenon. A scaling of the system does not change the equations much, which implies that a profile of a global solution can be found in a small neighborhood of a blowup point.
Since we require one parameter to tend to $4\pi$, it rules out the possibilities of fully bubbling profile for the whole system, which means around each bubbling disk, only the profile of a global solution for ONE equation occurs. One major difficulty in analysis is when multiple bubbling disks all tend to one blowup point. Here we invoke the following definition
 of $(\sigma_1(p),\sigma_2(p))$ in \cite{lwz-apde}: suppose $p$ is a blowup point,
$$\sigma_i(p)=\frac 1{2\pi}\lim_{\delta\to 0}\lim_{k\to \infty} \int_{B(p,\delta)}\rho_i^kh_ie^{u_1^k},\quad i=1,2. $$
Note that the limit of $k\to \infty$ is taken first, the second layer limit: $\lim_{\delta\to 0}$ is of secondary importance.
For regular $SU(3)$ Toda system all the possible formations are
$$(2,0),(0,2),(2,4),(4,2),(4,4). $$
Since we require $\rho_1^k\to 4\pi$, if $p$ is a blowup point of $u_1^k$, the only possible formations are $(2,0)$ or $(2,4)$. It has been proved in  \cite{musso} that it is possible to
construct a sequence of bubbling solutions $u^k=(u_1^k,u_2^k)$ such that $u_1^k$ and $u_2^k$ have a common blowup point $q^k$ and around $q^k$, the following spherical Harnack inequality holds for $u_2^k$:
\begin{equation}\label{har-u2}
u_2^k(x)+2\log |x-q^k|\le C,\quad \mbox{in}\quad B(q^k,\tau).
\end{equation}
 and if the system is scaled according to the maximum of $u_2^k$ around $q^k$, the scaled functions of $u_2^k$ converge to
$$\Delta v_2+2\lim_{k\to \infty} h_2^k(q^k)e^{v_2}=4\pi \delta_0,\quad \mbox{in}\quad \mathbb R^2, \quad \int_{\mathbb R^2}\lim_{k\to \infty}h_2^k(q^k)e^{v_2}<\infty. $$
On the other hand for the blowup type $(2,4)$, Ao-Wang \cite{ao} constructed $u^k=(u_1^k,u_2^k)$ such that $u_1^k$ and $u_2^k$ have a common blowup point, $q^k$ is a local maximum of $u_1^k$ and $u_2^k$ has two local blowup points $p_1^k,p_2^k$ such that
\begin{equation}\label{local-sym}
\lim_{k\to \infty} \frac{p_1^k-p_2^k}{\delta_k}=2\lim_{k\to \infty}(p_1^k-q^k)/\delta_k, \quad \delta_k=|p_1^k-q^k|\to 0.
\end{equation}
and
\begin{equation}\label{sp-har-2}
u_2^k(p_i^k)+2\log |x-q^k|\to \infty, \quad i=1,2.
\end{equation}

The main contribution of this article is to completely describe the profile of bubbling solutions when $(\rho_1^k,\rho_2^k)\to (4\pi,4\pi N)$. Our first conclusion is about possible formations of bubbling profiles, which says the two cases of bubbling collisions constructed in \cite{musso,ao} are the only two possible cases:

\begin{thm}\label{thm1} Let $u^k=(u_1^k,u_2^k)$ be a sequence of solutions of (\ref{main-2}) that satisfies (\ref{both-b}), (\ref{normal}) and (\ref{h1h2}). If $(\rho_1^k,\rho_2^k)\to (4\pi,4N\pi)$ for some $N\in \mathbb N$ and $N\ge 2$,
there are only three possible formations of bubbling solutions:
\begin{enumerate}
\item
$q^k\to q,\quad p_l^k\to p_l, \mbox{for } l=1,...,N$
where $q,p_1,...,p_N$ are $N+1$ distinct points, $q^k$ is a local maximum of $u_1^k$, $p_l^k$s are local maximums of $u_2^k$.
\item
$u_1^k$ has one blowup point $q$, $u_2^k$ has $N-1$ blowup points: $q,p_3,..,p_N$. For $\tau>0$ small, (\ref{har-u2}) holds for $u_2^k$.
\item
$u_1^k$ has one blowup point $q$, $u_2^k$ has $N-1$ blowup points $q,p_3,...,p_N$. In a small neighborhood of $q$ that excludes other blowup points, there are two local maximums of $u_2^k$ (denoted as $p_1^k$, $p_2^k$) and one local maximum of $u_1^k$ denoted as $q^k$, (\ref{local-sym}) and (\ref{sp-har-2}) hold for $p_1^k,p_2^k,q^k$.
\end{enumerate}
\end{thm}

We shall use cases 1,2,3 to describe the three alternatives in Theorem \ref{thm1}. Case one is an obvious possibility, case two and case three were constructed in \cite{musso} and \cite{ao} respectively.

Before stating the next result we fix some notations: we use $\lambda_{2,q}^k$ to denote the maximum of $u_2^k$ in $B(q^k,\tau)$ where $q=\lim_{k\to \infty}q^k$ is a common blowup point of $u_1^k$ and $u_2^k$, $q^k$ is a local maximum of $u_1^k$, $K(q^k)$ stands for the Gauss curvature at $q^k$. In case three we use $\delta_k$ to denote the distance between $p_1^k$ and $q^k$, finally we use $\nabla_1 G(x,y)$ to denote the differentiation with respect to the first component of $G(\cdot,\cdot)$.

For $\rho_2^k\to 4\pi$ or $8\pi$ we have the following results.

\begin{thm}\label{thm-simple} Let $u^k=(u_1^k,u_2^k)$ be a sequence of solutions of (\ref{main-2}) that satisfies (\ref{both-b}), (\ref{normal}) and (\ref{h1h2}). If $(\rho_1^k,\rho_2^k)\to (4\pi,4\pi)$, $u_1^k$ and $u_2^k$ each has one blowup point and these two points are distinct.
\end{thm}

In other words, only case one occurs for $(\rho_1^k,\rho_2^k)\to (4\pi,4\pi)$.

The following two theorems state simple situations for $(\rho_1^k,\rho_2^k)\to (4\pi,8\pi)$, the proofs of which are based on the main estimates for case two and case three in more general situations afterwards.

\begin{thm}\label{simple-2} Let $u^k=(u_1^k,u_2^k)$ be a sequence of solutions of (\ref{main-2}) that satisfies (\ref{both-b}), (\ref{normal}) and (\ref{h1h2}). Suppose $\rho_1^k\to 4\pi$ from below,
and $\rho_2^k\to 8\pi$. If $\nabla (\log h_1+2\log h_2)(q)+24\pi \nabla_1\beta(q,q)\neq 0$ on $M$, then $u_1^k$ has one blowup point, $u_2^k$ has two blowup points and these three points are all distinct.
\end{thm}
If we don't assume $\rho_1^k$ to tend to $4\pi$ from below, we have
\begin{thm}\label{simple-3}
Let $u^k=(u_1^k,u_2^k)$ be a sequence of solutions of (\ref{main-2}) that satisfies (\ref{both-b}), (\ref{normal}) and (\ref{h1h2}). Suppose $(\rho_1^k,\rho_2^k)\to (4\pi,8\pi)$ and
$\lim_{k\to \infty}\lambda_1^k/\lambda_2^k<4/3$. If $\nabla (\log h_1+2\log h_2)(q)+24\pi \nabla_1\beta(q,q)\neq 0$ on $M$, then $u_1^k$ has one blowup point, $u_2^k$ has two blowup points and these three points are all distinct.
\end{thm}

Our other main results include sharp estimates of the profile of bubbling solutions: Their pointwise estimate, location of blowup points, comparison of heights for each case. The readers are referred to each of the following sections for detailed statements. The main theorems in these sections
 provide crucial information about bubble interaction in terms of their heights, locations and proximity to threshold values. Just like in the single equation case, these results will play an important role in the construction of bubbles and the derivation of corresponding degree counting theorems. Even in the critical case
$(\rho_1,\rho_2)=(4\pi, 4\pi N)$, the results in this article are useful in the application of critical point at infinity. We will continue to carry out all the related theories in the future.

One major analytical difficulty related to Toda system is that maximum principles fail miserably. In this article we take advantage of the Harnack inequality proved in \cite{lwz-apde} and analyze the behavior of solutions using their spherical averages. At one point we need to approximate bubbling solutions accurately by global solutions. Since solutions to Toda systems may not have any symmetry, we use some key ideas in \cite{lwz-jems} to achieve this goal. In addition, many standard results for single Liouville equations will be used to obtain precise error estimates.

The organization of this article is as follows: In section two, all the estimates for case one are derived based on results for single equation. Then in section three we analyze the formation of bubbling coalition and then derive all the estimates for case three. Here the local maximums of $u_1^k$ and $u_2^k$ are proved to be on a line. As far as know this is the first time such an estimate has been proved for Toda system. In section four all the estimates for case two are established. One crucial ingredient in the proof of this section is a precise estimate for bubbling approximation, which depends heavily on the key idea in \cite{lwz-jems} for defeating the lack of symmetry and the non-degeneracy of the linearized Toda system established in \cite{lwy}. Finally we put the proofs of the theorems in the introduction in section five and explained the key idea in \cite{lwz-jems} for the simplest case in the appendix.

Notation: We use $B(p,r)$ to denote the ball with radius $r$ with center $p$. If $p$ is the origin we sometimes use $B_r$ instead of $B(0,r)$.

\section{Case one}

First we recall that case one means the only blowup point of $u_1^k$ is a not a blowup point of $u_2^k$. Thus there are $N+1$ blowup points all together. In this case we use the results of single Liouville equation to prove the following main result:
\begin{thm}\label{thm2}Let $u^k=(u_1^k,u_2^k)$ be a sequence of solutions of (\ref{main-2}) that satisfies (\ref{both-b}), (\ref{normal}) and (\ref{h1h2}). If $(\rho_1^k,\rho_2^k)\to (4\pi,4N\pi)$ and the situation described in case one of Theorem \ref{thm1} occurs, we have
$$
\rho_1^k-4\pi=2\frac{\Delta (\log h_1)(q^k)-2K(q^k)+4\pi N-8\pi}{h_1(q^k)}e^{-\lambda_1^k}\lambda_1^k \\
+O(e^{-\lambda_1^k})+O(e^{-\lambda_2^k}),
$$
where $q^k$ is a local maximum of $u_1^k$,
$$
\rho_2^k-4\pi N
=2\sum_{l=1}^N  e^{-\lambda_{2,l}^k}\lambda_{2,l}^k (\frac{\Delta (\log h_2)(p_l^k)-2K(p_l^k)+8N\pi-4\pi}{Nh_2(p_l^k)})
+O(e^{-\lambda_1^k})+O(e^{-\lambda_2^k}),
$$
where $p_l^k$( $l=1,...,N$) are local maximum points of $u_2^k$.
The locations of $q^k$, $p_1^k$,...,$p_N^k$ are related by
\begin{align*}
&\nabla (\log h_1)(q^k)-4\pi \sum_{l=1}^N \nabla_1 G(q^k,p_l^k)+8\pi \nabla_1\beta(q^k,q^k)\\
&=O(\lambda_1^ke^{-\lambda_1^k})+O(\lambda_2^ke^{-\lambda_2^k})\\
&\nabla (\log h_2)(p_l^k)-4\pi \nabla_1 G(p_l^k,q^k)+8\pi\sum_{s\neq l,s=1}^N\nabla_1 G(p_l^k,p_s^k)+8\pi \nabla_1\beta(p_l^k,p_l^k)\\
&=O(\lambda_1^ke^{-\lambda_1^k})+O(\lambda_2^ke^{-\lambda_2^k}),\quad l=1,...,N.
\end{align*}
 The comparison of the magnitudes of $u_2^k$ at each blowup is stated as
 \begin{align*}
\lambda_{2,l}^k+8\pi \beta(p_l^k,p_l^k)+8\pi \sum_{m\neq l}G(p_l^k,p_m^k)-12\pi G(p_l^k,q^k)+2\log h_2(p_l^k)
\nonumber \\
=\lambda_{2,s}^k+8\pi \beta(p_s^k,p_s^k)+8\pi \sum_{m\neq s}G(p_s^k,p_m^k)-12\pi G(p_s^k,q^k)+
2\log h_2(p_s^k)\nonumber\\
+O(\lambda_2^k\lambda_2^ke^{-\lambda_2^k})+O(\lambda_2^ke^{-\lambda_1^k}), \quad \mbox{ for } \quad l\neq s.
\end{align*}
\end{thm}
Here we note that in case one there is no direct relation between $\lambda_1^k$ and $\lambda_2^k$.  This phenomenon changes in case two and case three.

In order to prove those estimates in Theorem \ref{thm1} we need to focus on one equation frequently. We first apply the following standard procedure to single
out one equation near a blowup point: Suppose $q$ is a blowup point of $u_1^k$ and not a blowup point of $u_2^k$, we let
$$f_2^k(x)=-\int_M G(x,y)\rho_2^k h_2e^{u_2^k}dV_g(y) $$
be a solution of
$$\Delta_g f_2^k=\rho_2^k(h_2e^{u_2^k}-1),\quad \mbox{ in }\quad M. $$
Using $f_2^k$ we can remove $u_2^k$ from the equation for $u_1^k$:
\begin{equation}\label{eq-q-loc}
\Delta_g (u_1^k-f_2^k)+2\rho_1^k h_1 e^{u_1^k}=2\rho_1^k, \quad \mbox{ in }\quad M.
\end{equation}
Now we use $ds^2=e^{\phi}|dx|^2$ to denote the isothermal coordinate in a neighborhood of $q$: Suppose $q^k$ is a local maximum of $u_1^k$ ($\lim_{k\to \infty}q^k=q$), in local coordinates,
\begin{equation}\label{iso-c}
\phi(0)=0,\quad |\nabla \phi(0)|=0,\quad \Delta \phi(0)=-2K(q^k)
\end{equation}
where $K(q^k)$ is the Guass curvature at $q^k$. With the conformal covariant property of $\Delta_g$ in surfaces ( $\Delta_g =e^{-\phi}\Delta$)  we write (\ref{eq-q-loc}) in a neighborhood of $q^k$ as
$$\Delta (u_1^k-f_2^k)+2\rho_1^kh_1 e^{\phi} e^{u_1^k}=2\rho_1^k e^{\phi},\quad \mbox{ in }\quad B_{\tau} $$
where $\tau>0$ is chosen that $B(q^k,2\tau)$ does not contain other blowup points.

For applications later it is more convenient to make the equation homogeneous and make the solution having no oscillation on $\partial B_{\tau}$. For this purpose we set
$f_3^k$ be a solution of
\begin{equation}\label{f3k-e}
\Delta f_3^k(x)=e^{\phi},\quad \mbox{in}\quad B_{\tau},\quad f_3^k(0)=|\nabla f_3^k(0)|=0,
\end{equation}
 and let $\psi_1^k$ be a harmonic function on $B_{\tau}$ that eliminates the oscillation of $u_1^k-f_2^k-2\rho_2^kf_3^k$ on $\partial B_{\tau}$:
\begin{equation}\label{psi-1}
\left\{\begin{array}{ll}\Delta \psi_1^k=0,\quad \mbox{ in}\quad B_{\tau},\\
\psi_1^k(x)=\frac{1}{2\pi \tau}\int_{\partial B_{\tau}} (u_1^k-f_2^k-2\rho_2^kf_3^k)dS,\quad x\in \partial B_{\tau}.
\end{array}
\right.
\end{equation}
By the mean value property of harmonic functions $\psi_1^k(0)=0$. Later in application it is clear that $u_1^k$, $f_2^k$ and $f_3^k$ all have finite oscillation on $\partial B_{\tau}$ so all the derivatives of $\psi_1^k$ are bounded in $B_{\tau/2}$. Lumping $f_2^k$, $f_3^k$ and $\psi_1^k$ with $u_1^k$ we write (\ref{eq-q-loc}) as
\begin{equation}\label{til-u-1k}
\Delta \tilde u_1^k+2h_1^ke^{\tilde u_1^k}=0,\quad \mbox{in}\quad B_{\tau}
\end{equation}
where
\begin{equation}\label{t-u1-h}
\tilde u_1^k=u_1^k-f_2^k-2\rho_1^k f_3^k-\psi_1^k, \quad h_1^k=\rho_1^k h_1 e^{f_2^k+2\rho_2^kf_3^k+\psi_1^k},
\end{equation}
we here emphasize that $\tilde u_1^k$ is a constant (depending on $k$) on $\partial B_{\tau}$.

Before deriving sophisticated estimate we start with a crude one:
\begin{lem}\label{rough-rho12}
\begin{equation}\label{rough-2}
\rho_1^k-4\pi=O(\lambda_1^ke^{-\lambda_1^k}),\quad \rho_2^k-4\pi N=O(\lambda_2^ke^{-\lambda_2^k}).
\end{equation}
\end{lem}

\noindent{\bf Proof of Lemma \ref{rough-rho12}:}

Here we recall that $q^k$ is a local maximum of $u_1^k$ ($\lambda_1^k=u_1^k(q^k)=\max_M u_1^k$), which is $0$ in the local coordinate system.
Let $\tilde q^k$ be the maximum point of $\tilde u_1^k$, here we cite a standard result for single equation \cite{gluck, cl1, zhangcmp}, $\tilde q^k=O(e^{-\lambda_1^k})$,
$\tilde u_1^k(\tilde q^k)-\lambda_1^k=O(e^{-\lambda_1^k})$,
\begin{equation}\label{good-ap-2}
\nabla (\log  h_1^k)(0)=O(\lambda_1^ke^{-\lambda_1^k}),
\end{equation}
\begin{equation}\label{good-ap-1}
\tilde u_1^k(x)=\log \frac{e^{ \lambda_1^k}}{(1+e^{ \lambda_1^k}\frac{2 h_1^k(0)}8|x-\tilde q^k|^2)^2}+O( \lambda_1^k\lambda_1^ke^{-\lambda_1^k}),\quad x\in B_{\tau}.
\end{equation}
Using (\ref{good-ap-1}) and (\ref{good-ap-2}) in the following computation, we have
$$2\int_{B_{\tau}} h_1^ke^{\tilde u_1^k}=8\pi +O(\lambda_1^ke^{-\lambda_1^k}). $$

Since
$$2\rho_1^k\int_{B(q^k,\tau)}h_1e^{u_1^k}dV_g=2\int_{B_{\tau}} h_1^ke^{\tilde u_1^k}$$ and
$u_1^k=-\lambda_1^k+O(1)$ for $x\in M\setminus B(q^k,\tau)$, we see that the estimate of $\rho_1^k-4\pi$ in (\ref{rough-2}) holds. The estimate of $\rho_2^k-4\pi N$ can be derived in a similar fashion since all the blowup points of $u_1^k$ and $u_2^k$ are distinct. The influence of $u_1^k$ can be removed from the equation for $u_2^k$ around a blowup point and everything can be carried out as before.
Lemma \ref{rough-rho12} is established. $\Box$

\medskip

In order to derive more precise estimates we use the expansion of bubbling solution for Liouville equation, as employed in the proof of Lemma \ref{rough-rho12}, to have
\begin{equation}\label{f2k-q}
f_2^k(q^k)=-4\pi \sum_{l=1}^NG(q^k,p_l^k)+O(\lambda_2^ke^{-\lambda_2^k}),
\end{equation}
and
\begin{equation}\label{f2k-qd}
\nabla f_2^k(q^k)=-4\pi \sum_{l=1}^N \nabla_1 G(q^k,p_l^k)+O(\lambda_2^ke^{-\lambda_2^k}).
\end{equation}
Here we note that if $N>2$, $u_2^k$ has $N-2$ blowup points other that $q^k$ in case two and case three. Let $p_{2,l}^k$ ($l=3,...,N$) be these blowup points, by results for single equation
$$\lambda_{2,l}^k=\lambda_{2,s}^k+O(1),\quad l,s=3,...,N $$
where $\lambda_{2,l}^k$ is $u_2^k(p_{2,l}^k)$, $l=3,..,N$.

Since certain errors occur very frequently we use the following notations:
$$SE_1=O(\lambda_1^ke^{-\lambda_1^k})+O(\lambda_2^ke^{-\lambda_2^k}),\quad SE_2=O(e^{-\lambda_1^k})+O(e^{-\lambda_2^k}).$$
Taking advantage of these new notations in direct computation, we have, from the definition of $h_i^k$ in (\ref{t-u1-h}), that
\begin{eqnarray}\label{t-u-h}
h_1^k(0)=4\pi h_1(q^k)e^{-4\pi\sum_{l=1}^NG(q^k,p_l^k)}+SE_1,\\
\tilde u_1^k(0)=u_1^k(q^k)+4\pi \sum_{l} G(q^k,p_l^k)+SE_1. \nonumber
\end{eqnarray}

Now we compute
$$\int_{B(0,\tau)} h_1^ke^{\tilde u_1^k}=\int_{\Omega_k} h_1^k(\epsilon_ky)e^{v_k}dy$$
where $$v_k(y)=\tilde u_1^k(\epsilon_ky)+2\log \epsilon_k,$$
$\epsilon_k=e^{-\frac 12\tilde u_1^k(0)}$,
$\Omega_k=B(0,\epsilon_k^{-1}\tau)$,
the leading term in the expansion of $v_k$ is
$$U_k(y)=\log \frac{1}{(1+\frac{ h_1^k(0)}4|y|^2)^2},  $$
and the expansion is
$$v_k(y)=U_k(y)+O(\epsilon_k^2)(\log (2+|y|))^2. $$

Using the expansion in the calculation we have
\begin{align}\label{tem-3}
&\int_{B(0,\tau)} h_1^ke^{\tilde u_1^k}
=\int_{B(0,\tau \epsilon_k^{-1})}( h_1^k)(\epsilon_k y)e^{v_k}dy\\
=&4\pi+O(\epsilon_k^2)+\frac 12\sum_{i=1}^2 \partial_{ii}( h_1^k)(0)\epsilon_k^2\int_{B(0,\tau \epsilon_k^{-1})}y_i^2e^{U}dy \nonumber\\
=&4\pi +\frac{8\pi}{ h_1^k(0)^2}\Delta h_1^k(0)\epsilon_k^2\log \epsilon_k^{-1}+O(\epsilon_k^2).\nonumber
\end{align}
Note that the term involving $\nabla h_1^k(0)$ is part of an error because of symmetry.

The Pohozaev identity for single equation  gives
$$\nabla (\log h_1^k)(0)=O(e^{-\lambda_1^k}\lambda_1^k), $$
which immediately implies
$$\Delta h_1^k(0)=h_1^k(0)\Delta (\log h_1^k)(0)+SE_1. $$
Thus (\ref{tem-3}) can be written as
\begin{equation}\label{tem-4}
\int_{B(0,\tau)} h_1^ke^{\tilde u_1^k}-4\pi=8\pi\frac{\Delta(\log h_1^k)(0)}{h_1^k(0)}\epsilon_k^2\log \frac{1}{\epsilon_k}+SE_1.
\end{equation}

Going back to the definition of $h_1^k$ we have

\begin{eqnarray*}
\Delta (\log h_1^k)(0)=\Delta (\log h_1)(q^k)+\Delta \phi+\Delta f_2^k(q^k)+2\rho_1^k\Delta f_3^k(0)\\
=\Delta (\log h_1)(q^k)-2K(q^k)-4\pi N+8\pi+SE_1.
\end{eqnarray*}

Let $\tilde \rho_1^k=\rho_1^k\int_{B(q^k,\tau)}h_1e^{u_1^k}$, using
$$\epsilon_k^2=e^{-\tilde u_1^k(0)}
=e^{-u_1^k(q^k)-4\pi\sum_lG(q^k,p_l^k)+O(\lambda_2^ke^{-\lambda_2^k})},$$ and the corresponding expression for $\log 1/\epsilon_k$ in (\ref{tem-4}), we have

\begin{equation*}
\tilde \rho_1^k-4\pi
=2\frac{\Delta (\log h_1)(q^k)-2K(q^k)+4\pi N-8\pi}{h_1(q^k)}e^{-\lambda_1^k}\lambda_1^k+SE_2.
\end{equation*}

Since $\int_{M\setminus B(q^k,\tau)}\rho_1^kh_1 e^{u_1^k}dV_g=O(e^{-\lambda_1^k})$, we have
\begin{equation}\label{rho-1-case-1}
\rho_1^k-4\pi=2\frac{\Delta (\log h_1)(q^k)-2K(q^k)+4\pi N-8\pi}{h_1(q^k)}e^{-\lambda_1^k}\lambda_1^k +SE_2.
\end{equation}

The following lemma determines the location of $q^k$:
\begin{lem}\label{lem-qk}
\begin{equation}\label{case-1-loc}
\nabla (\log h_1)(q^k)-4\pi \sum_{l=1}^N \nabla_1 G(q^k,p_l^k)+8\pi \nabla_1\beta(q^k,q^k)=SE_1.
\end{equation}
\end{lem}

\noindent{\bf Proof of Lemma \ref{lem-qk}:}

Let $\xi\in \mathbb S^1$, the following Pohozaev identity for $\tilde u_1^k$ holds:
\begin{equation}\label{pi-tu1}
\int_{B_{\tau}}\partial_{\xi} h_1^k e^{\tilde u_1^k}=\int_{\partial B_{\tau}} e^{\tilde u_1^k}  h_1^k (\xi \cdot \nu)
+\partial_{\nu} \tilde u_1^k \partial_{\xi} \tilde u_1^k-\frac 12 |\nabla \tilde u_1^k|^2(\xi \cdot \nu).
\end{equation}
First we observe that the first term on the right hand side is minor because its order is $SE_2$. In order to evaluate other terms on the right hand side of (\ref{pi-tu1}) we consider the Green's representation of $u_1^k-f_2^k$ around $q^k$ (see (\ref{pi-tu1}) :

$$u_1^k(x)-f_2^k(x)=\bar u_1^k-\bar f_2^k+\int_M G(x,\eta)2\rho_1^k h_1e^{f_2^k}e^{u_1^k-f_2^k}dV, $$
we use the expansion of $u_1^k-f_2^k$ to obtain, for $x$ away from $q^k$, that
$$u_1^k(x)-f_2^k(x)=\overline{ u_1^k-f_2^k}+ 8\pi  G(x,q^k)+SE_1. $$
In view of the fact that $\bar f_2^k=0$ based on its definition, we have
$$u_1^k-f_2^k-2\rho_1^k f_3^k=\bar u_1^k-4\log |x|+8\pi \beta(x,q^k)-2\rho_1^k f_3^k +SE_1, \quad x\in B_{\tau}\setminus \{0\}. $$
The right hand side without $SE_1$
$$\tilde G(x)=-4\log |x|+8\pi \beta(x,q^k)-8\pi f_3 $$
 is a harmonic function in $B_{\tau}\setminus \{0\}$, which means
 $$\tilde u_1^k=u_1^k-f_2^k-2\rho_1^k f_3^k-\psi_1^k=\bar u_1^k+C+SE_1, \quad |x|=\tau, $$
 and
 $$\nabla \tilde u_1^k=SE_1,\quad |x|=\tau. $$
 Using this in the Pohozaev identity we have

\begin{equation}\label{q-case-1a}
\nabla (\log  h_1^k)(0)+8\pi \nabla_1\beta(q^k,q^k)-4\pi \sum_{l}\nabla_1 G(q^k,p_l^k)=SE_1.
\end{equation}
By the definition of $h_1^k$ we see that (\ref{case-1-loc}) holds. Lemma \ref{lem-qk} is established. $\Box$

\medskip

Now we determine the leading terms of $\rho_2^k-4\pi N$. Define $\rho_{2,l}^k$ by the integration of $\rho_2 h_2^ke^{u_2^k}$ over $B(p_l^k,\tau)$ for some small $\tau>0$ (so that $B(p_1^k,2\tau)$ does not contain other blowup points). Let
$$f_1^k(x)=-\int_MG(x,\eta) \rho_1^k h_1 e^{u_1^k}dV_g(\eta). $$
Then for $x$ away from $q^k$ we have
\begin{equation}\label{f1-k-p}
f_1^k(x)=-4\pi G(x,q^k)+O(\lambda_1^ke^{-\lambda_1^k})
\end{equation}
and
\begin{equation}\label{f1-k-d}
\nabla f_1^k(x)=-4\pi \nabla_1 G(x,q^k)+O(\lambda_1^ke^{-\lambda_1^k}).
\end{equation}
The equation for $f_1^k$
$$\Delta_g f_1^k=\rho_1^k(h_1 e^{u_1^k}-1), $$
gives rise to
\begin{equation}\label{u2-f1}
\Delta_g(u_2^k-f_1^k)+2\rho_2^k(h_2 e^{f_1^k}e^{u_2^k-f_1^k}-1)=0.
\end{equation}
Later we shall use (\ref{u2-f1}) to determine the relations between $\lambda_{2,l}^k=u_2^k(p_l^k)$.
Now we write (\ref{u2-f1}) as

$$\Delta_g(u_2^k-f_1^k)+2\rho_2^k h_2 e^{u_2^k}=2\rho_2^k. $$
In isothermal coordinates around $p_l^k$, we have

$$\Delta (u_2^k-f_1^k)+2\rho_2^k e^{\phi} h_2 e^{f_1^k}e^{u_2^k-f_1^k}=2\rho_2^k e^{\phi} $$
for $\phi(0)=|\nabla \phi(0)|=0$, $\Delta \phi(0)=-2K(p_l^k)$.

Let $f_3$ be a solution of
$$\Delta f_3=e^{\phi}, \quad f_3(p_l^k)=|\nabla f_3(p_l^k)|=0 $$
and set
$$\tilde u_2^k=u_2^k-f_1^k-2\rho_2^k f_3^k-\psi_2^k, \quad  h_2^k=\rho_2^k h_2 e^{\phi+f_1^k+2\rho_2^kf_3^k+\psi_2^k}, $$
where $\psi_2^k$ is a harmonic function determined by the average of $u_2^k-f_1^k-2\rho_2^k f_3^k$ on $\partial B_{\tau}$.
The equation for $\tilde u_2^k$ is
$$\Delta \tilde u_2^k+ 2h_2^k e^{\tilde u_2^k}=0. $$
The values of $\tilde u_2^k$ and $h_2^k$ at the origin are
\begin{equation}\label{u2pl}
\tilde u_2^k(0)=\lambda_{2,l}^k+4\pi G(p_l^k,q^k)+O(\lambda_1^ke^{-\lambda_1^k})
\end{equation}
and
\begin{equation}\label{t-h20}
 h_2^k(0)=4\pi Nh_2(p_l^k)e^{-4\pi G(p_l^k,q^k)}+SE_1.
\end{equation}

Green's representation of $u_2^k$ gives
$$(u_2^k-f_1^k)(x)=\bar u_2^k-\bar f_1^k+\int_M G(x,\eta)2\rho_2^k h_2 e^{u_2^k}dV_{\eta}. $$

After using results for single equation we have, for $x$ away from $p_l^k$,
$$u_2^k(x)-f_1^k(x)=\bar u_2^k+\sum_{s=1}^N8\pi G(x, p_s^k)+SE_1. $$

Using the Pohozaev identity for single equation, we obtain the location of $p_l^k$ as follows:
\begin{eqnarray}
\nabla (\log h_2)(p_l^k)-4\pi \nabla_1 G(p_l^k,q^k)+8\pi\sum_{s\neq l,s=1}^N\nabla_1 G(p_l^k,p_s^k)+8\pi \nabla_1\beta(p_l^k,p_l^k) \nonumber \\
=SE_1,\quad l=1,...,N.\label{case-1-loc-2}
\end{eqnarray}

Let $\rho_{2,l}^k=\int_{B(p_l^k,\tau)}\rho_2^k h_2 e^{u_2^k}$.
Then we have
\begin{equation*}
\rho_{2,l}^k=\frac 12 \int_{B(p_l^k,\tau)}2\rho_2^k h_2 e^{u_2^k}=\frac 12 \int_{B_{\tau}} 2 h_2^ke^{\tilde u_2^k}.
\end{equation*}

Let
$$\epsilon_k=e^{-\frac 12 \tilde u_2^k(0)},\quad \epsilon_k^2=e^{-\lambda_{2,l}^k-4\pi G(p_l^k,q^k)+O(\lambda_1^ke^{-\lambda_1^k})}, $$
and
$$v_k(y)=\tilde u_2^k(\epsilon_k y)+2\log \epsilon_k, $$
Then the leading term in the expansion of $v_k$ is
$$\log \frac{1}{(1+\frac{ h_2^k(0)}4|y|^2)^2}. $$
Then following the same computation as before we have
\begin{eqnarray*}
\rho_{2,l}^k=\frac 12 (8\pi+\frac 12\epsilon_k^2\Delta  h_2^k(0)\int_{B(0,\epsilon_k^{-1}\tau)}|y|^2e^{v(y)}dy+O(\epsilon_k^2)\\
=4\pi +8\pi \epsilon_k^2\log \frac{1}{\epsilon_k}\Delta  h_2^k (0)/ h_2^k(0)^2+O(\epsilon_k^2).
\end{eqnarray*}

By the definition of $h_2^k$ we have
\begin{equation}\label{del-t-h}
\Delta (\log h_2^k)(0)=\Delta (\log h_2)(p_l^k)-2K(p_l^k)+8\pi N-4\pi+O(\lambda_1^ke^{-\lambda_1^k}).
\end{equation}
The Pohozaev identity for single equation gives
$$\nabla (\log h_2^k)(0)=O(\lambda_2^ke^{-\lambda_2^k}), $$
which leads to
$$\Delta (\log h_2^k)(0)=\frac{\Delta h_2^k(0)}{ h_2^k(0)}+O(\lambda_2^ke^{-\lambda_2^k}). $$

Using the definition of $\tilde u_2^k(p_l^k)$ and $\epsilon_k$, we have
\begin{equation*}
\rho_{2,l}^k-4\pi
=2 e^{-\lambda_{2,l}^k}\lambda_{2,l}^k(\frac{\Delta (\log h_2)(p_l^k)-2K(p_l^k)+8N\pi-4\pi}{Nh_2(p_l^k)})+SE_2.
\end{equation*}

Thus
\begin{equation}\label{case-1-rho-2}
\rho_2^k-4\pi N
=2\sum_{l=1}^N  e^{-\lambda_{2,l}^k}\lambda_{2,l}^k (\frac{\Delta (\log h_2)(p_l^k)-2K(p_l^k)+8N\pi-4\pi}{Nh_2(p_l^k)})+SE_2.
\end{equation}

Now we use the results for single equation to derive the mutual relationship between $\lambda_{2,l}^k$ and $\lambda_{2,s}^k$($l\neq s$). First it is easy to see that
the only difference between any two of them is $O(1)$.

From
$$u_2^k(x)=\bar u_2^k+\int_M G(x,\eta)(2\rho_2^kh_2 e^{u_2^k}-\rho_1^kh_1 e^{u_1^k})dV_{\eta}, $$
we focus on $x=p_l^k$ and have
\begin{align*}
&\lambda_{2,l}^k=u_2^k(p_l^k)=\bar u_2^k+\int_{B(p_l^k,\delta)}G(p_l^k,\eta)(2\rho_2^kh_2 e^{u_2^k}-\rho_1^kh_1 e^{u_1^k})dV_{\eta}\\
&+\sum_{s\neq l}\int_{B(p_s^k,\tau)}+\int_{B(q^k,\tau)}+\int_{M\setminus (\cup_l B(p_l^k,\tau)\cup B(q^k,\tau))}\\
&=\bar u_2^k-\frac 1{2\pi}\int_{B(p_l^k,\tau)}\log |\eta | 2\rho_2^kh_2 e^{u_2^k}+8\pi \beta(p_l^k,p_l^k)+\sum_{s\neq l}8\pi G(p_l^k,p_s^k)
-4\pi G(p_l^k,q^k)\\
&+SE_1.
\end{align*}
To evaluate the second term on the right hand side, we use
$$\epsilon_k=e^{-\frac 12 \tilde u_2^k(0)},\quad v_k(y)=\tilde u_2^k(\epsilon_ky)+2\log \epsilon_k, $$
where $\tilde u_2^k(0)$ satisfies (\ref{u2pl}).

\begin{eqnarray*}
(-\frac 1{2\pi})\int_{B(p_l^k,\tau)}\log |\eta | 2\rho_2^k h_2e^{u_2^k}=
(-\frac 1{2\pi })\int_{\Omega_k}(\log |\eta_1|+\log \epsilon_k)2 h_2^k(\epsilon_k \eta_1)e^{v_k(\eta_1)}d\eta_1\\
=-4\log \epsilon_k-\frac 1{2\pi}\int_{\Omega_k}\log |\eta_1| 2h_2^k(\epsilon_k \eta_1)e^{v_k(\eta_1)}d\eta+O(\lambda_2^k\lambda_2^ke^{-\lambda_2^k})+O(\lambda_2^ke^{-\lambda_1^k})
\end{eqnarray*}
where $\Omega_k=B(0,\tau \epsilon_k^{-1})$, $h_2^k$ is as before, the scaled version of $\tilde u_2^k$ is $v_k$, which satisfies
$$v_k(y)=\log \frac{1}{(1+\frac{h_2^k(0)}4 |y|^2)^2}+O((\log \epsilon_k)^2). $$
It is elementary to verify that
\begin{eqnarray*}
(-\frac 1{2\pi})\int_{\Omega_k}\log |\eta_1 | 2 h_2^k(0)e^{v_k}
=2\log \frac{h_2^k(0)}4+O(\lambda_2^ke^{\lambda_2^k})\\
=2\log (\pi N h_2(p_2^k))-8\pi G(p_l^k,q^k)+SE_1.
\end{eqnarray*}

Note that the following equality is used in the evaluation above:
$$\int_0^{\infty}\log r\frac{r}{(1+r^2)^2}dr=0.$$

Putting all the information together we have
\begin{eqnarray*}
0=\bar u_2^k+\lambda_{2,l}^k+8\pi \beta(p_l^k,p_l^k)+\sum_{s\neq l}8\pi G(p_l^k,p_s^k)-4\pi G(p_l^k,q^k)\\
+2\log N\pi+2\log h_2(p_l^k)+O(\lambda_2^k\lambda_2^k e^{-\lambda_2^k})+O(\lambda_2^ke^{-\lambda_1^k}).
\end{eqnarray*}Here

Thus for $l\neq s$, equating $\bar u_k$ leads to
\begin{align}
\lambda_{2,l}^k+8\pi \beta(p_l^k,p_l^k)+8\pi \sum_{m\neq l}G(p_l^k,p_m^k)-4\pi G(p_l^k,q^k)+2\log h_2(p_l^k)
\nonumber \\
=\lambda_{2,s}^k+8\pi \beta(p_s^k,p_s^k)+8\pi \sum_{m\neq s}G(p_s^k,p_m^k)-4\pi G(p_s^k,q^k)+
2\log h_2(p_s^k)\nonumber\\
+O(\lambda_2^k\lambda_2^ke^{-\lambda_2^k})+O(\lambda_2^ke^{-\lambda_1^k}). \label{comp-c-1}
\end{align}

\medskip

\noindent{\bf Proof of Theorem \ref{thm2} }.   The difference of $\rho_1^k$ and $4\pi$ can be found in (\ref{rho-1-case-1}), $\rho_2^k-4N\pi$ is stated in (\ref{case-1-rho-2}). The location that $q^k$ satisfies, as a local maximum of $u^k$ is derived in (\ref{q-case-1a}) and for each local maximum $p_l^k$ of $u_2^k$, the location is determined in (\ref{case-1-loc-2}).
Finally the comparison of the magnitudes of local maximum values of $u_2^k$ is established in (\ref{comp-c-1}).
Theorem \ref{thm2} is established. $\Box$

\section{Case three}

In this section we establish sharp estimates for case three stated in Theorem \ref{thm1}. First we recall that $q^k$ is the only blowup point of $u_1^k$ and $q=\lim_{k\to \infty}q^k$ is also a blowup point of $u_2^k$. Moreover, $u_2^k$ has two local maximum points $p_1^k,p_2^k$, both tending to $q$ but the spherical Harnack inequality around $q^k$ is violated:
$$u_2^k(p_i^k)+2\log |p_i^k-q^k|\to \infty, \quad i=1,2. $$
We use $\delta_k=|p_1^k-q^k|$. Our main results for case three in the statement of Theorem \ref{thm1} are as follows. First for $N=2$, we have
\begin{thm}\label{thm4} In case three and $N=2$, the following conclusions hold:
\begin{enumerate}
\item $\lambda_1^k+(4+d_k)\log \delta_k\to \infty$,  $\lambda_2^k+2\log \delta_k\to \infty$, where $d_k=O(e^{-\lambda_2^k}\delta_k^{-2})$
\item $\rho_1^k-4\pi=(D+o(1))e^{-\lambda_1^k}\delta_k^{-4-d_k}+O(\lambda_1^ke^{-\lambda_1^k })$ where $D>0$ is a constant.
\item $\rho_2^k-8\pi=O(\lambda_2^ke^{-\lambda_2^k})+O(e^{-\lambda_2^k}\delta_k^{-2})$.
\item
\begin{align*}
&\nabla \log h_1(q^k)+2\nabla \log h_2(q^k)+24\pi\nabla_1\beta(q^k,q^k)\\
=&O(\lambda_1^ke^{-\lambda_1^k}\delta_k^{-3})+O(e^{-\lambda_1^k}\delta_k^{-4-d_{k}})+O(\lambda_2^ke^{-\lambda_2^k}\delta_k^{-3})+O(\delta_k).
\end{align*}
\end{enumerate}
\end{thm}

In case three and $N\ge 3$, we have
\begin{thm}\label{thm-case-3-N}
For case three and $N\ge 3$, the following conclusions hold: For some $d_k=O(e^{-\lambda_{2,q}}\delta_k^{-2})$ we have
\begin{enumerate}
\item $\lambda_1^k+(4+d_k)\log \delta_k\to \infty$,
\item $\lambda_2^k+(4+2d_k)\log \delta_k\to \infty$,
\item $\lambda_{2,q}^k-(2+2d_k)\log \delta_k=\lambda_2^k+O(1)$,
\item $\rho_1^k-4\pi=(D+o(1)) e^{-\lambda_1^k}\delta_k^{-4-d_k}\big (\int_M h_1 exp(-4\pi\sum_{l=3}^NG(x,p_l^k))dV_g \big)\, +O(\lambda_1^ke^{-\lambda_1^k })$ for some $D>0$.
\item $\rho_2^k-4\pi N=O(e^{-\lambda_{2,q}^k})+O(\lambda_{2,q}^ke^{-\lambda_{2,q}^k}\delta_k^{-2})$.
\item The locations of $q^k$, $p_l^k$ ($l\ge 3$) are determined by
\begin{align*}
\nabla (\log h_1+2\log h_2)(q^k)+24\pi \nabla_1\beta(q^k,q^k)+12\pi\sum_{l\ge 3}\nabla_1G(q^k,p_l^k)\\
=O(\lambda_1^k e^{-\lambda_1^k}\delta_k^{-3})+O(e^{-\lambda_1^k}\delta_k^{-4-d_k})+O(\lambda_{2,q}^ke^{-\lambda_{2,q}^k}\delta_k^{-3})+O(\delta_k),
\end{align*}
and
\begin{align*}
\nabla \log h_2(p_l^k)+12\pi\nabla_1 G(p_l^k,q^k)+8\pi \nabla_1\beta(p_l^k,p_l^k)+8\pi \sum_{s\neq l}\nabla_1 G(p_l^k,p_s^k)\\
=O(e^{-\lambda_1^k}\delta_k^{-4-d_k})+O(\lambda_{2,q}^ke^{-\lambda_{2,q}^k})+O(e^{-\lambda_{2,q}^k}\delta_k^{-2}),\quad l=3,..,N.
\end{align*}
Finally the comparison of heights can be found in
\begin{align*}
\lambda_{2,q}^k-2\log \delta_k+\log \frac{h_2(q^k)}8+12\pi \beta(q^k,q^k)+8\pi \sum_{l\ge 3}G(q^k,p_l^k) \nonumber \\
=\lambda_{2,l}^k+\log \frac{h_2(p_l^k)}4+12\pi G(p_l^k,q^k)+8\pi \beta (p_l^k,p_l^k)+\sum_{s\neq l,s\ge 3}8\pi G(p_l^k,p_s^k)\nonumber \\
+O(\delta_k)+O(\lambda_2^k\lambda_2^ke^{-\lambda_2^k}\delta_k^{-4-2d_k})+O(\lambda_1^ke^{-\lambda_1^k}\delta_k^{-4-d_k}\log \delta_k^{-1}).
\end{align*}
\end{enumerate}
\end{thm}

Note that the comparison
between $\lambda_{2,l}^k$ and $\lambda_{2,m}^k$ for $l\neq m$ ($l,m\in \{3,...,N\}$) can be obtained by the last estimate in Theorem \ref{thm-case-3-N}. Also the $d_k$ in Theorem \ref{thm-case-3-N} is $o(1)$ because of (2) and (3).

Case three is about  bubbling disks colliding into one point. Since $\rho_1^k\to 4\pi$, at least around the bubbling disk that contains the maximum of $u_1^k$, we have
$$\int_{B(q^k,\tau)}2\rho_1^kh_1e^{u_1^k}\ge 8\pi +o(1). $$
Using the fact that $\int_{B(q^k,\tau)}h_1e^{u_1^k}\le 1$,
we see immediately that there is only one bubbling disk for $u_1^k$. we recall the definition of $\sigma_i$ in \cite{lwz-apde}:
$$\sigma_i=\frac 1{2\pi}\lim_{\delta\to 0}\lim_{k\to \infty} \int_{B(q^k,\delta)}h_ie^{u_i^k}dV_g,\quad i=1,2. $$
Since the classification of $(\sigma_1,\sigma_2)$ for $SU(3)$ Toda system is
$$(2,0),(0,2),(2,4),(4,2),(4,4),$$
the only type that can have colliding bubbling disks is $(2,4)$, which has two possible formations: (1) Either one bubbling disk
that contains a partial blowup profile of $u_1^k$ with two other bubbling disks of partial blowup of $u_2^k$, or, (2) after scaling according to the magnitude of $u_2^k$, the scaled function of $u_2^k$ tending to
$$\Delta v+2e^v=4\pi \delta_0,\quad \mbox{ in }\quad \mathbb R^2, \quad \int_{\mathbb R^2}e^v<\infty. $$

 We first consider the case that $(2,4)$ consists of three bubbling disks of
 $$(2,0),(0,2),(0,2)$$
 respectively. Let $q^k$ and $p_1^k$ and $p_2^k$ be the centers of these three bubbling disks. We first prove that $q^k$ is roughly the mid-point of $p_1^k$ and $p_2^k$. Recall that $\lambda_{2,q}^k$ is the maximum $u_2^k$ in $B(q^k,\tau)$.

 Since $q^k$ is the common blowup point for $u_1^k$ and $u_2^k$, we now consider the local version of the system in a neighborhood of $q^k$: First using the $\phi$ function corresponding to the isothermal coordinates at $q^k$ we write the system as
$$\left\{\begin{array}{ll}
\Delta u_1^k+2\rho_1^kh_1 e^{\phi}e^{u_1^k}-\rho_2^k h_2e^{\phi}e^{u_2^k}=(2\rho_1^k-\rho_2^k)e^{\phi},\\
\Delta u_2^k-\rho_1^kh_1 e^{\phi}e^{u_1^k}+2\rho_2^k h_2e^{\phi}e^{u_2^k}=(-\rho_1^k+2\rho_2^k)e^{\phi}, \quad \mbox{in}\quad B_{\tau}.
\end{array}
\right.
$$
Note that $\phi$ satisfies
\begin{equation}\label{phi-cur}
\phi(0)=|\nabla \phi(0)|=0,\quad \Delta \phi(0)=-2K(q^k).
\end{equation}
 Then we use the following two functions to eliminate the right hand side:
 \begin{align}\label{case3-f}
\Delta f_1^k=(2\rho_1^k-\rho_2^k)e^{\phi},\quad f_1^k(0)=|\nabla f_1^k(0)|=0,\\
\Delta f_2^k=(-\rho_1^k+2\rho_2^k)e^{\phi},\quad f_2^k(0)=|\nabla f_2^k(0)|=0.\nonumber
\end{align}
Let $\psi_i^k$ be the harmonic function that eliminates the oscillation of $u_i^k-f_i^k$ on $B(q^k,\tau)$ for some $\tau>0$ small. By setting
\begin{equation}\label{case3-uihi}
\tilde u_i^k=u_i^k-f_i^k-\psi_i^k, \quad  h_i^k=\rho_i^kh_ie^{\phi}e^{f_i^k+\psi_i^k},\quad i=1,2
\end{equation}
we write the system in $B(0,\tau)$ as
\begin{equation}\label{sys-q}
\left\{\begin{array}{ll}
\Delta \tilde u_1^k+2 h_1^k(x)e^{\tilde u_1^k}- h_2^ke^{\tilde u_2^k}=0,\\
\Delta \tilde u_2^k- h_1^k(x)e^{\tilde u_1^k}+2 h_2^k e^{\tilde u_2^k}=0,\quad \mbox{in}\quad B_{\tau},
\end{array}
\right.
\end{equation}
and the definition of $h_i^k$ gives
\begin{eqnarray}\label{loc-q-12}
\nabla (\log h_1^k)(0)=\nabla (\log h_1)(q^k)+\nabla \psi_1^k(0)\nonumber \\
\nabla (\log h_2^k)(0)=\nabla (\log h_2)(q^k)+\nabla \psi_2^k(0)\nonumber \\
\Delta (\log h_1^k)(0)=\Delta (\log h_1)(q^k)+8\pi-4N\pi-2K(q^k)+o(1)\\
\Delta (\log h_2^k)(0)=\Delta (\log h_2)(q^k)+8\pi N-4\pi-2K(q^k)+o(1)\nonumber
\end{eqnarray}
where $o(1)=O(|\rho_1^k-4\pi|)+O(|\rho_2^k-4\pi N|)$.

Since there is only bubbling disk for $u_1^k$, the maximum of $u_1^k$ is in $B(q^k,\tau)$:
$\lambda_1^k=\max_{B(q^k,\tau)}u_1^k$. Note that the maximum of $\tilde u_1^k$ in $B_{\tau}$ is only $O(e^{-\lambda_1^k})$ different from $\lambda_1^k$ (see \cite{gluck, zhangcmp}).

\medskip

Before we deduce more specific locations of $Q_i$, we recall that the concept of group, defined in \cite{lwz-apde}, describes a few bubbling disks of comparable distances to one another but far away to bubbling disks not in the group.
\begin{lem} \label{mid-pt} Suppose $\delta_k=|q^k-p_1^k|\le |q^k-p_2^k|$. Note that since they are in a group, $|q^k-p_2^k|\le C|q^k-p_1^k|$.
Let
$$v_i^k(y)=\tilde u_i^k(q^k+\delta_k y)+2\log \delta_k,\quad i=1,2 $$
and $Q_i^k$ be the images of $p_i^k$ after scaling ($i=1,2$, $Q_i^k\to Q_i$), then we have $Q_2=-Q_1$.
\end{lem}

\noindent{\bf Proof of Lemma \ref{mid-pt}:}

The system after scaling becomes
$$\left\{\begin{array}{ll}
\Delta v_1^k(y)+2h_1^k(\delta_k y)e^{v_1^k(y)}-h_2^k(\delta_k y)e^{v_2^k(y)}=0,\quad |y|\le R_k,\\
\Delta v_2^k(y)-h_1^k(\delta_k y)e^{v_1^k(y)}+2h_2^k(\delta_k y)e^{v_2^k(y)}=0,\quad |y|\le R_k.
\end{array}
\right.
$$
for some $R_k\to \infty$.
Let $Q_1^k$ and $Q_2^k$ be the images of $p_1^k$ and $p_2^k$ respectively, the definition of $Q_1^k$ gives $|Q_1^k|=1$. For $Q_2^k$ we know $|Q_2^k|\sim 1$. For $\delta>0$ small,  $B_{\delta}$, $B(Q_1^k,\delta)$ and $B(Q_2^k,\delta)$ are disjoint and in these disks we can see the profile of bubbling solutions of a single Liouville equation.  If we let $\mu_1^k$, $\mu_2^k$ and $\bar \mu_2^k$ be the maximum values of $v_1^k$ in $B(0,\delta)$, magnitude of $v_2^k$ in $B(Q_1^k,\delta)$ and $B(Q_2^k,\delta)$ respectively. It is obvious that $\mu_2^k=\bar \mu_2^k+O(1)$. Using Green's representation it is also clear that $v_i^k$ has bounded oscillation around each of the bubbling disks. If we consider the location of blowup at the origin, it satisfies
\begin{equation}\label{tem-0}
\delta_k \nabla (\log h_1^k)(0)+\nabla \phi_1^k(0)=O(\mu_1^ke^{-\mu_1^k}),
\end{equation}
where $\phi_1^k$ is the harmonic function determined by the value of $v_1^k$ on $\partial B(0,1/2)$:
\begin{equation}\label{tem-1}
\phi_1^k=-2\log |y-Q_1^k|-2\log |y-Q_2^k|+C+o(1).
\end{equation}
Then
\begin{equation}\label{tem-2}
\nabla \phi_1^k(0)=\frac{2Q_1}{|Q_1|^2}+2\frac{Q_2}{|Q_2|^2}+o(1).
\end{equation}

Combining (\ref{tem-0}) and (\ref{tem-2}) we obtain $Q_2^k=-Q_1^k+o(1)$.
$\Box$

\medskip

With the new information revealed in Lemma \ref{mid-pt} we now consider the equation for $(v_1^k,v_2^k)$ on $B(0,\tau\delta_k^{-1})$:
$$
\left\{\begin{array}{ll}\Delta v_1^k+2h_1^k(\delta_ky)e^{v_1^k}-h_2^k(\delta_ky)e^{v_2^k}=0,\\
\Delta v_2^k(y)-h_1^k(\delta_ky)e^{v_1^k}+2h_2^k(\delta_ky)e^{v_2^k}=0.\quad \mbox{in}\quad B(0,\delta_k^{-1}).
\end{array}
\right.
$$

 Let
$$\mu_2^k=\lambda_{2,q}^k+2\log \delta_k\to \infty, $$
then the standard uniform estimate for single equation (\cite{li-cmp}) gives a rough estimate of $v_2^k$ around $Q_1$:,
$$v_2^k(y)=-\mu_2^k+O(1), \quad y\in \partial B(Q_1,\delta). $$
For more detailed analysis of $v_i^k$ we need to consider the equation of their spherical averages $\bar v_i^k(r)$:
\begin{equation}\label{bar-v2}
\frac{d}{dr}\bar v_2^k(r)=-\frac{2\sigma_2^k(r)-\sigma_1^k(r)}r,\quad \bar v_2^k(2)=-\mu_2^k+O(1).
\end{equation}
$$v_2^k(y)=-\mu_2^k-(6+o(1))\log |y|+O(1),\quad 2<|y|<\delta_k^{-1}, $$
where
$$\sigma_i^k(r)=\frac 1{2\pi}\int_{B_r}h_i^k(\delta_ky)e^{v_i^k(y)}dy,\quad i=1,2. $$
Here we note that
\begin{equation}\label{sm-u1}
\int_{M\setminus B(q^k,2\delta_k)} h_1e^{u_1^k} dV_g=o(1).
\end{equation}
The reason is by the result of the single Liouville equation we already have
$$2\int_{B_{1/2}}h_1^k(\delta_k y)e^{v_1^k(y)}=8\pi+o(1), $$
this fact, together with $\rho_1^k\to 4\pi$ and $\int_M h_1 e^{u_1^k}=1$, proves (\ref{sm-u1}).

Let
$$d_{1,k}=\sigma_1^k(2)-2,\quad d_k=\sigma_2^k(2)-4, $$
equation (\ref{sm-u1}) further implies that
\begin{equation}\label{mu-1k-big}
\mu_1^k+(2+d_k)\log \delta_k\to \infty, \quad d_k=O(e^{-\mu_2^k}).
\end{equation}
To see this  we consider the equation of $\bar v_1^k(r)$. Since $v_2^k$ has very fast decay, it is easy to see from $\sigma_1^k(r)\le 2+o(1)$ (because $\rho_1^k\to 4\pi$) that
$$\bar v_2^k(r)=-\frac{2\sigma_2^k(r)-\sigma_1^k(r)}r\le -(6+o(1))/r,\quad 2<r<\tau \delta_k^{-1}. $$
We further observe from the results for single equation that
\begin{equation}\label{soft-1}
\sigma_1^k(2)=2+O(\mu_1^k e^{-\lambda_1^k})+O(e^{-\mu_1^k})
\end{equation}
and
\begin{equation}\label{soft-2}
 \sigma_2^k(2)=4+O(\mu_2^k e^{-\lambda_{2,q}^k})+O(e^{-\mu_2^k}).
\end{equation}
The reason for (\ref{soft-1}) is
$$\sigma_1^k(2)=2+C\Delta h_1^k(0) \delta_k^2e^{-\mu_1^k}+O(e^{-\mu_1^k}) $$
where $\delta_k^2$ comes from the scaling. Thus (\ref{soft-1}) holds. (\ref{soft-2}) can be proved similarly.
Lemma 2.1 of \cite{lwz-apde} says that $v_i^k$ is $O(1)$ different from its spherical average away from bubbling area. Thus we consider the equation of $\bar v_1^k(r)$ first. Recall that
$$\frac{d}{dr}\bar v_1^k(r)=-\frac{2\sigma_1^k(r)-\sigma_2^k(r)}r=-\frac{2(\sigma_1^k(r)-2)-(\sigma_2^k(r)-4)}r. $$
Let
$$\epsilon_k=\sigma_2^k(2)-2\sigma_1^k(2)=(\sigma_2^k(2)-4)-2(\sigma_1^k(2)-2)=d_k-2d_{k,1},$$
then (\ref{soft-1}) and (\ref{soft-2}) yield
\begin{equation}\label{ep-est}
\epsilon_k=O(\mu_1^ke^{-\lambda_1^k})+O(e^{-\mu_1^k})+O(\mu_2^ke^{-\lambda_{2,q}^k})+O(e^{-\mu_2^k}),
\end{equation}
Let
$$\tilde \sigma_i^k(r)=\frac 1{2\pi}\int_{B_r\setminus B_2}h_i^k(\delta_ky)e^{v_i^k(y)}dy=\sigma_i^k(r)-\sigma_i^k(2),\quad i=1,2. $$
Then the equation for $\bar v_1^k(r)$ can be written as
$$\frac{d}{dr}\bar v_1^k(r)=\frac{\epsilon_k-2\tilde \sigma_1^k(r)+\tilde \sigma_2^k(r)}r,\quad 2<r<\tau \delta_k^{-1}. $$

Because of the decay of $v_2^k$ we have
\begin{equation}\label{t-sig-2-small}
\tilde \sigma_2^k(r)=O(e^{-\mu_2^k}r^{-4+o(1)}),\quad 2<r<\tau \delta_k^{-1}.
\end{equation}
Now for $r>2$,
$$\frac{d}{dr}\bar v_1^k(r)=\frac{\epsilon_k-2\tilde \sigma_1^k(r)+\tilde \sigma_2^k(r)}r,\quad 2<r<\tau \delta_k^{-1} $$
with $\bar v_1^k(2)=-\mu_1^k+O(1)$.
Integrating the equation for $\bar v_1^k(r)$ above, we have, using Harnack inequality,
$$\bar v_1^k(r)\le -\mu_1^k+\epsilon_k \log r+C, $$
which further gives, by Lemma 2.1 of \cite{lwz-apde},
$$e^{v_1^k(y)}\le Ce^{-\mu_1^k}r^{\epsilon_k}, \quad |y|=r\in (2, \tau \delta_k^{-1}). $$
This estimate of $v_1^k$ gives a crude estimate of $\tilde \sigma_1^k(r)$:
\begin{equation}\label{better-s-1}
\tilde \sigma_1^k(r)\le Ce^{-\mu_1^k}r^{2+\epsilon_k},
\end{equation}
which leads to a lower bound of $\bar v_1^k(r)$ as
$$\frac{d}{dr}\bar v_1^k(r)\ge \frac{\epsilon_k}r-2Ce^{-\mu_1^k}r^{1+\epsilon_k}. $$
Consequently
$$v_1^k(y)\ge -\mu_1^k+\epsilon_k \log r-C_1e^{-\mu_1^k}r^{2+\epsilon_k}-C_2,\quad |y|=r\in (2, \tau \delta_k^{-1}), $$
which is equivalent to
$$e^{v_1^k(y)}\ge C_1 e^{-\mu_1^k}r^{\epsilon_k}e^{-C_2e^{-\mu_1^k}r^{2+\epsilon_k}}, \quad |y|=r. $$
Thus we obtain a lower bound of $\tilde \sigma_1^k(r)$:
\begin{align*}
\tilde \sigma_1^k(r)\ge Ce^{-\mu_1^k}\int_2^rs^{1+\epsilon_k}e^{-Ce^{-\mu_1^k}s^{2+\epsilon_k}}ds\\
=C_3\int_{C_2 e^{-\mu_1^k}}^{C_2 e^{-\mu_1^k}r^{2+\epsilon_k}}e^{-t}dt
\ge C e^{-\mu_1^k}r^{2+\epsilon_k}.
\end{align*}
In particular for $r\sim \delta_k^{-1}$ we have
\begin{equation}\label{about-s-1}
\mu_1^k+(2+\epsilon_k)\log \delta_k\to \infty,
\end{equation}
because otherwise $\rho_1^k$ would not tend to $4\pi$. As an immediate consequence $\mu_1^k+(2+o(1))\log \delta_k\to \infty$, which means the two terms related to $\mu_1^k$ in (\ref{ep-est}) are completely harmless and the estimate of $\tilde \sigma_1^k(r)$ becomes
\begin{equation}\label{sig-new-out}
\tilde \sigma_1^k(r)\sim e^{-\mu_1^k}r^{2+d_k}.
\end{equation}
Then it is easy to see that
\begin{equation}\label{v-1k-outside}
v_1^k(y)=-\mu_1^k+d_k\log r+O(1),\quad 2<r<\tau \delta_k^{-1}.
\end{equation}
In particular
 (\ref{mu-1k-big}) holds.

One quick observation based on (\ref{mu-1k-big}) is that we now have a more precise estimate of $v_2^k$ away from the bubbling area:
\begin{equation}\label{better-v2k}
v_2^k(y)=-\mu_2^k-(6+2d_k)\log |y|+O(1),\quad 2<|y|<\tau \delta_k^{-1}.
\end{equation}

\medskip

Now we prove more precise description of $p_1^k$ and $p_2^k$, the two local maximum points of $v_2^k$ around $q^k$.
From
$$\Delta v_2^k+2h_2^k(\delta_ky)e^{v_2^k}=h_1^k(\delta_ky)e^{v_1^k(y)} $$
we write
$$v_2^k(y)=\int_{\Omega_k}G_k(y,\eta)(2h_2^ke^{v_2^k}-h_1^ke^{v_1^k(\eta)})d\eta+v_2^k|_{\partial \Omega_k} $$
where $\Omega_k$ is the re-scaled domain
$\Omega_k=B(0,\tau \delta_k^{-1})$ for some $\tau>0$ small, and $v_2^k$ is constant on $\partial \Omega_k$, $G_k$ is the Green's function with respect to Dirichlet boundary condition:
\begin{equation}\label{green-k}
G_k(y,\eta)=-\frac 1{2\pi}\log |y-\eta |+H(y,\eta)
\end{equation}
with
$$H(y,\eta)=\frac 1{2\pi}\log \frac{|\eta |}{\tau \delta_k^{-1}}|\frac{\tau^{2}\delta_k^{-2}\eta}{|\eta |^2}-y |. $$
For $y\in B_3\setminus (B(Q_1^k,1/4)\cup B(Q_2^k,1/4)\cup B(0,1/4))$ we have, by estimates for single equations,
\begin{align}\label{case-3-tem-1}
\nabla v_2^k(y)&=\int_{\Omega_k}\nabla_1 G_k(y,\eta)(2h_2^ke^{v_2^k}-h_1^ke^{v_1^k})d\eta \nonumber \\
&=-4\frac{y-Q_2^k}{|y-Q_2^k|^2}-4\frac{y-Q_1^k}{|y-Q_1^k|^2}+2\frac{y}{|y|^2}+E
\end{align}
where
$$E=O(\mu_2^ke^{-\mu_2^k})+O(\mu_1^ke^{-\mu_1^k})+O(\delta_k^2).$$
From (\ref{case-3-tem-1}) it is easy to see that the harmonic function that eliminates the oscillation of $v_2^k$ around $Q_1^k$ is
$$\phi_{2,1}^k(y)=-4\log |y-Q_2^k|+2\log |y|+C+E, $$
which satisfies
\begin{equation}\label{osi-1}
\nabla \phi_{2,1}^k(Q_1^k)=(-4)\frac{Q_1^k-Q_2^k}{|Q_1^k-Q_2^k|^2}+2\frac{Q_1^k}{|Q_1^k|^2}+E.
\end{equation}
Similarly, if $\phi_{2,2}^k$ is the harmonic function that eliminates the oscillation of $v_2^k$ around $Q_2^k$, we have
\begin{equation}\label{osi-2}
\nabla \phi_{2,2}^k(Q_2^k)=(-4)\frac{Q_2^k-Q_1^k}{|Q_2^k-Q_1^k|^2}+2\frac{Q_2^k}{|Q_2^k|^2}+E.
\end{equation}
Using this in the Pohozaev identity around $Q_1^k$, we have
\begin{equation}\label{loc-pi-q}
\delta_k \nabla (\log h_2^k)(\delta_k Q_1^k)+(-4)\frac{Q_1^k-Q_2^k}{|Q_1^k-Q_2^k|^2}+2\frac{Q_1^k}{|Q_1^k|^2}=E.
\end{equation}
Since $Q_1^k\in \partial B_1$, we use $Q_1^k=e^{i\theta_k}$. So $Q_2^k=(-\cos \theta_k,-\sin \theta_k)+o(1)$. We set it as
$$Q_2^k=(-\cos \theta_k+a,-\sin \theta_k+b). $$
Using this expression of $Q_2^k$ in (\ref{loc-pi-q}), we observe by elementary computation that
\begin{align*}
&\frac{Q_1^k-Q_2^k}{|Q_1^k-Q_2^k|^2}\\
=&\frac{(2\cos \theta_k-a, 2\sin \theta_k-b)}{4-4a\cos \theta_k -4b\sin\theta_k +O(a^2+b^2)}\\
=&(\frac 12 \cos \theta_k+\frac{a}{4}\cos 2\theta_k+\frac{b}{4}\sin 2\theta_k, \frac 12 \sin \theta_k+\frac{a}{4}\sin 2\theta_k-\frac{b}4 \cos 2\theta_k)
+O(a^2+b^2).
\end{align*}
Now (\ref{loc-pi-q}) reads
\begin{align}\label{det-a-1}
\delta_k \partial_1 (\log h_2^k)(0)-a \cos 2\theta_k-b\sin 2\theta_k+O(a^2+b^2)=E \nonumber \\
\delta_k \partial_2 (\log h_2^k)(0)-a \sin 2\theta_k+b\cos 2\theta_k+O(a^2+b^2)=E.
\end{align}
Using the fact that $a,b=o(1)$ we first observe that
\begin{equation}\label{ab-pre}
a,b=O(\delta_k)+E.
\end{equation}
Inserting (\ref{ab-pre}) in (\ref{det-a-1}) we have
\begin{align}\label{a-and-b}
a=\delta_k\partial_1(\log h_2^k)(0)\cos 2\theta_k+\delta_k \partial_2(\log h_2^k)(0)\sin 2\theta_k+E\\
b=-\delta_k\partial_2(\log h_2^k)(0)\cos 2\theta_k+\delta_k \partial_1(\log h_2^k)(0)\sin 2\theta_k+E, \nonumber
\end{align}

Now we consider the equation for $v_1^k$:
$$\Delta v_1^k+2 h_1^k(\delta_ky)e^{v_1^k}= h_2^k(\delta_k y)e^{v_2^k}. $$
Let
$$f_{2d}^k(y)=\frac{1}{2\pi}\int_{\Omega_k} \log |y-\eta | h_2^ke^{v_2^k}d\eta. $$
After computation we have
$$f_{2d}^k(y)=2\log |y-Q_2^k|+2\log |y-Q_1^k|+E $$
for $y$ around $0$, and
\begin{align*}
\nabla f_{2d}^k(0)&=-\frac{2Q_2^k}{|Q_2^k|^2}-2\frac{Q_1^k}{|Q_1^k|^2}+E\\
&=(2a cos 2\theta_k+2b \sin 2\theta_k,-2b\cos 2\theta_k+2a \sin 2\theta_k)+E\\
&=2\delta_k(\partial_1 (\log h_2^k)(0),\partial_2 (\log h_2^k)(0))+E.
\end{align*}
Pohozaev identity for single equation gives
$$\delta_k\nabla (\log h_1^k)(0)+\nabla f_{2,d}^k(0)=E, $$
which is
\begin{equation}\label{loc-case-3}
\nabla \log h_1^k(0)+2\nabla \log h_2^k(0)=O(\delta_k^{-1}E).
\end{equation}

Now we discuss the rough behavior of $u_2^k$ around $q^k$. The asymptotic behavior of $v_2^k$ is (\ref{better-v2k}).
For $|y|\sim \delta_k^{-1}$, which is $|x-q^k|\sim 1$ or $|y|\sim \delta_k^{-1}$,
$$\tilde u_2^k(x)+2\log \delta_k=v_2^k(y)=-\mu_2^k+(6+2d_k)\log \delta_k+O(1). $$
Replacing $\mu_2^k$ by $\lambda_{2,q}^k$, we have
$$u_2^k(x)=-\lambda_{2,q}^k+(2+2d_k)\log \delta_k+O(1),\quad x\in \partial B(q^k,\tau). $$
If $N\ge 3$ and $p_l^k$ is another blowup point for $u_2^k$ with local maximum value $\lambda_{2,l}^k$, we know from the results for single equations that around $p_l^k$,
$$\tilde u_2^k(x)=-\lambda_{2,l}^k+O(1),\quad |x-p_l^k|=\tau. $$
Thus the comparison of the value of $u_2^k$ away from bubbling disks gives
\begin{equation}\label{lambda-2-q}
-\lambda_{2,q}^k+(2+2d_k)\log \delta_k=-\lambda_{2,l}^k+O(1),\quad l=3,...,N.
\end{equation}

 In order to determine the harmonic functions that cancel the oscillation of $u_i^k$ on $\partial B(q^k,\tau)$, we use the Green's representation of $u_i^k$ around $q^k$:
$$u_1^k(x)=\bar u_1^k+\int_M G(x,y)(2\rho_1^kh_1e^{u_1^k}-\rho_2^k h_2 e^{u_2^k})dV_g(y) $$

For $x\in M\setminus B(q^k,\tau)$, if $N=2$, it is easy to obtain
\begin{equation}\label{u1k-2}
u_1^k(x)=\bar u_1^k+O(\mu_2^ke^{-\lambda_{2,q}^k})+O(e^{-\mu_2^k})+O(\mu_1^ke^{-\lambda_1^k})+O(e^{-\lambda_1^k}\delta_k^{-d_k-4}).
\end{equation}
Indeed,
we evaluate the integral into integration over $B(q^k,2\delta_k)$ and $M\setminus B(q^k,2\delta_k)$. The integration of $G(x,\eta) \rho_2^kh_2e^{u_2^k(\eta)}$ over these two regions gives
$$8\pi G(x,q^k)+O(e^{-\mu_2^k})+O(\mu_2^ke^{-\lambda_{2,q}^k})$$ as the sum.  The evaluation of
$$\int_{B(q^k,2\delta_k)} G(x,\eta)2\rho_1^k h_1 e^{u_1^k} dV_g $$ gives
$8\pi G(x,q^k)+O(e^{-\lambda_1^k}\mu_1^k)+O(e^{-\mu_1^k})$ and the integration of the same integrand over $M\setminus B(q^k,2\delta_k)$ is $O(e^{-\mu_1^k}\delta_k^{-2-d_k})$. Thus (\ref{u1k-2}) holds.
If $N\ge 3$, we have
\begin{align}\label{u1k-N}
u_1^k(x)&=\bar u_1^k-4\pi\sum_{l\ge 3} G(x,p_l^k)\nonumber\\
&+O(\mu_2^ke^{-\lambda_{2,q}^k})+O(e^{-\mu_2^k})+O(\mu_1^ke^{-\lambda_1^k})+O(e^{-\mu_1^k}\delta_k^{-2-d_k})+O(\lambda_2^ke^{-\lambda_2^k}).
\end{align}
For simplicity we set
$$SE_3=O(e^{-\lambda_1^k}\delta_k^{-4-d_k})+O(\mu_2^ke^{-\lambda_{2,q}^k})+O(e^{-\mu_{2}^k})+O(\mu_1^ke^{-\lambda_1^k}).$$

By studying the ODE of $\bar v_1^k(r)$ and the Green's representation of $u_1^k$, we obtain easily that
$$\bar u_1^k=-\mu_1^k+(2+d_k)\log \delta_k^{-1}+C+SE_3=-\lambda_1^k-(4+d_k)\log \delta_k+C+SE_3, $$
which means, for $D=e^C$,
$$e^{\bar u_1^k}=(D+o(1))e^{-\lambda_1^k}\delta_k^{-4-d_k}, \quad N=2 $$
and
$$e^{u_1^k}=(D+o(1))e^{-\lambda_1^k}\delta_k^{-4-d_k}e^{-4\pi \sum_l G(x,p_l^k)},\quad N\ge 3. $$

With the estimate of $e^{\bar u_1^k}$ we evaluate $\rho_1^k-4\pi$ as follows: First for a small neighborhood of $q^k$, using results for single equation we obtain
\begin{equation}\label{tem-11}
\int_{B(q^k,\delta_k/2)}\rho_1^kh_1e^{u_1^k}dV_g=4\pi+O(\mu_1^ke^{-\mu_1^k}\delta_k^2)+O(e^{-\mu_1^k}).
\end{equation}
Note that $\delta_k^2$ comes from differentiating coefficient function twice. Then the integration over $M\setminus B(q^k,N\delta_k)$ gives, using Green's representation formula of $u_1^k$,
\begin{equation}\label{tem-12}
\int_{M\setminus B(q^k,N\delta_k)}\rho_1^kh_1e^{u_1^k}dV_g=(D+o(1))e^{-\mu_1^k}\delta_k^{-2-d_k}, \quad \mbox{ for } \quad N=2
\end{equation}
and
\begin{equation}\label{tem-13}
\int_{M\setminus B(q^k,N\delta_k)}\rho_1^kh_1e^{u_1^k}dV_g=(D+o(1))e^{-\mu_1^k}\delta_k^{-2-d_k}\int_M h_1(x) e^{-4\pi \sum_l G(x,p_l^k)}dV_g(x)
\end{equation}
for $N\ge 3$. Note that $-4\pi G(x,p_l^k)$ has a good sign near $p_l^k$ and is integrable near $p_l^k$.
Thus $\rho_1^k-4\pi$ in Theorem \ref{thm4} and Theorem \ref{thm-case-3-N} can be derived from (\ref{tem-11}), (\ref{tem-12}) and (\ref{tem-13}).

To evaluate $\rho_2^k-4N\pi$, we first consider the integration around $q^k$: Let $\tilde \rho_2^k=\int_{B(q^k,2\delta_k)}\rho_2^kh_2e^{u_2^k}dV_g$, we first claim that
\begin{equation}\label{rho-2-q}
\tilde \rho_2^k-8\pi=O(\delta_k^{-2}e^{-\lambda_{2,q}^k})+O(\lambda_{2,q}^ke^{-\lambda_{2,q}^k}).
\end{equation}
Indeed,
$$\int_{B(p_1^k,\delta_k/2)}\rho_2^kh_2e^{u_2^k}dV_g=\int_{B(Q_1^k,\frac 12)}h_2^ke^{v_2^k}=4\pi+O(\mu_2^ke^{-\mu_2^k}\delta_k^2)+O(e^{-\mu_2^k}). $$
Here we note that $\delta_k^{2}$ comes from the Laplacian of the coefficient function. The integration of over $B(p_2^k,\delta_k/2)$ can be evaluated similarly.
Since
$$v_2^k(y)=-\mu_2^k+O(1) \quad \mbox{on}\quad \partial B_2$$
and $v_2^k$ has a very fast decay, it is easy to see that
$$\int_{B(q^k,\tau/2)\setminus B(q^k,2\delta_k)}\rho_2^kh_2e^{u_2^k}dV_g=\int_{B(0,\tau \delta_k^{-1}/2)\setminus B_2}h_2^k(\delta_ky)e^{v_2^k(y)}dy=O(e^{-\mu_2^k}). $$
Thus (\ref{rho-2-q}) follows.
Next Green's representation for $u_2^k$ gives
\begin{align}\label{u2-e}
&u_2^k(x)=\bar u_2^k+\int_M G(x,y)(2\rho_2^kh_2 e^{u_2^k}-\rho_1^k h_1e^{u_1^k})dV_{y}\\
&=\bar u_2^k+ 12 \pi G(x,q^k)+\sum_{l=3}^N 8\pi G(x, p_l^k)+SE_3,\quad x\in M\setminus B(q^k,\tau/2).\nonumber
\end{align}

Thus, from (\ref{u1k-2}) and (\ref{u1k-N}) we see that the derivatives of the harmonic function that eliminates the oscillation of $u_1^k$ around $q^k$ satisfy
\begin{equation}\label{loc-q-34}
\nabla \psi_1^k(q^k)=\left\{\begin{array}{ll}SE_3,\quad \mbox{for } N=2,\\
-4\pi \sum_{l\ge 3}\nabla_1G(q^k,p_l^k)+SE_3,\quad \mbox{for}\quad N\ge 3,
\end{array}
\right.
\end{equation}
Note that $\nabla f_i^k(0)=0$ is used.
Similarly from (\ref{u2-e}) we have
\begin{equation}\label{2nd-har}
\nabla \psi_2^k(q^k)=\left\{\begin{array}{ll} 12\pi \nabla_1 \beta (q^k,q^k)+SE_3,\quad N=2,\\
12\pi \nabla_1\beta(q^k,q^k)+8\pi \sum_{l\ge 3}\nabla_1 G(q^k,p_l^k)+SE_3 \quad N\ge 3,
\end{array}
\right. \nonumber
\end{equation}

Combining (\ref{loc-q-34}) with (\ref{loc-case-3}) we have, for $N=2$, that
\begin{equation}\label{loc-q-2}
\nabla (\log h_1+2\log h_2)(q^k)+24\pi \nabla_1\beta(q^k,q^k)=SE_3+E \delta_k^{-1}.
\end{equation}
The vanishing estimate for $N\ge 3$ becomes
\begin{equation}\label{loc-q-f}
\nabla (\log h_1+2\log h_2)(q^k)+24\pi \nabla_1\beta(q^k,q^k)+12\pi\sum_{l\ge 3}\nabla_1G(q^k,p_l^k)
=SE_3+E \delta_k^{-1}.
\end{equation}
Thus the location statements in Theorem \ref{thm4} and Theorem \ref{thm-case-3-N} are justified.


Note that the spreading of bubbling circles implies $\lambda_{2,q}^k+2\log \delta_k\to \infty$, which gives
$$\lambda_2^k+(4+2d_k)\log \delta_k\to \infty ,\quad \mbox{if}\quad N\ge 3 $$
based on (\ref{lambda-2-q}).
In order to determine the relation between $\lambda_{2,q}^k$ with $\lambda_{2,l}^k$, we use Green's representation of $u_2^k$. It is easy to see that
\begin{equation}\label{u2k-qk}
u_2^k(x)=\bar u_2^k-6\log |x-q^k|+12\pi \beta(x,q^k)+8\pi \sum_{l\ge 3} G(x,p_l^k)+SE_3 \log \delta_k^{-1}
\end{equation}
for $x\in B(q^k,\tau)\setminus B(q^k,\tau/2)$. To determine the relation between $\bar u_2^k$ and $\lambda_{2,q}^k$, we consider the expansion of $v_2^k$ around $Q_2^k$ and without loss of generality we assume that $v_2^k(Q_2^k)=\mu_2^k$. Then for
 $\frac 18<|y-Q_2^k|<1/2$, standard results for single equation give
\begin{align}\label{v2k-e1}
v_2^k(y)-\psi_{2,v}^k&=-\mu_2^k-2\log \frac{h_2(q^k)}4-4\log |y-Q_2^k|+O(\delta_k)\\
&+O(\mu_2^k\mu_2^ke^{-\mu_2^k})+O(\mu_1^ke^{-\mu_1^k}).\nonumber
\end{align}
where $\psi_{2,v}^k$ is a harmonic function that eliminates the oscillation of $v_2^k$ on $\partial B(Q_2^k,1/2)$ and satisfies $\psi_{2,v}^k(Q_2^k)=0$.
For simplicity we set
$$SE_4=O(\delta_k)+O(\mu_2^k\mu_2^ke^{-\mu_2^k})+O(\mu_1^ke^{-\mu_1^k}).$$

Now let $\Omega_k=B(0,\tau \delta_k^{-1})$ and we consider the expression of $v_2^k$ on $\Omega_k$ by the Green's representation formula:
\begin{equation}\label{v2k-e2}
v_2^k(y)=\int_{\Omega_k}G_k(y,\eta)(2h_2^ke^{v_2^k}-h_1^ke^{v_1^k})+v_2^k|_{\partial \Omega_k}.
\end{equation}
Using $G_k(y,\eta)$ as in (\ref{green-k}), we have, for $y\not \in B(Q_1^k,1/2)\cup B(Q_2^k,1/2)$
\begin{align}\label{v2k-qk}
v_2^k(y)&=-4\log |y-Q_1^k|-4\log |y-Q_2^k|+2\log |y|+6\log (\tau \delta_k^{-1})\\
&+v_2^k|_{\partial B(0,\tau \delta_k^{-1})}+SE_3\log \delta_k^{-1}, \nonumber
\end{align}
which determines
$$\psi_{2,v}^k=-4\log |y-Q_1^k|+2\log |y|+4\log 2+SE_3\log \delta_k^{-1}+SE_4 $$
and
\begin{align}\label{v2k-e8}
v_2^k-\psi_{2,v}^k=-4\log |y-Q_2^k|+6\log \tau-6\log \delta_k +v_2^k|_{\partial \Omega_k}-4\log 2\nonumber\\
+SE_3\log \delta_k^{-1}+SE_4.
\end{align}

The comparison of $v_2^k$ in (\ref{v2k-e1}) and (\ref{v2k-e8}) gives the value of $v_2^k$ on $\partial \Omega_k$:
\begin{align}\label{v2k-bry}
v_2^k|_{\partial \Omega_k}=-\mu_2^k-6\log \tau+6\log \delta_k-2\log \frac{h_2(q^k)}8+SE_3 \log (\delta_k^{-1})+SE_4\\
=-\lambda_{2,q}^k-6\log \tau+4\log \delta_k-2\log \frac{h_2(q^k)}8+SE_3 \log \delta_k^{-1}+SE_4. \nonumber
\end{align}
Then we immediately obtain
\begin{equation}\label{bar-u2-tau}
\tilde u_2^k|_{\partial B(q^k,\tau)}=-\lambda_{2,q}^k-6\log \tau+2\log \delta_k-2\log \frac{h_2(q^k)}8+SE_3\log \delta_k^{-1}+SE_4.
\end{equation}
Recall that
$$\tilde u_2^k=u_2^k-f_2^k-\psi_2^k. $$
Therefore
\begin{align}\label{u2-tau}
u_2^k|_{\partial B(q^k,\tau)}=f_2^k+\psi_2^k-\lambda_{2,q}^k-6\log \tau+2\log \delta_k-2\log \frac{h_2(q^k)}8 \nonumber \\
+SE_3\log \delta_k^{-1}+SE_4.
\end{align}
The equation of $f_2^k+\psi_2^k$ determines that
\begin{align*}
f_2^k(x)+\psi_2^k(x)=12\pi \beta (x, q^k)-12\pi \beta(q^k,q^k)+8\pi(\sum_{l=3}^N G(x,p_l^k)-G(q^k,p_l^k))\\
+SE_3\log \delta_k^{-1}+SE_4.
\end{align*}
Evaluating the Green's representation of $u_2^k$ near $q^k$ gives, for $|x-q^k|\sim \tau$,
\begin{align}\label{bar-u2-tau2}
&u_2^k(x)=\bar u_2^k+\int_MG(x,\eta)(2\rho_2^kh_2e^{u_2^k}-\rho_1^kh_1e^{u_1^k})dV_g(\eta)\\
&=\bar u_2^k+12\pi G(x,q^k)+\sum_{l=3}^N 8\pi G(x,p_l^k)+SE_3\log \delta_k^{-1}+SE_4\nonumber\\
&=\bar u_2^k+f_2^k+\psi_2^k-6\log |x-q^k|+12\pi \beta(q^k,q^k)
+8\pi\sum_{l=3}^N G(q^k,p_l^k)+SE_3\log \delta_k^{-1}+SE_4. \nonumber
\end{align}
where the expression of $f_2^k+\psi_2^k$ is used.  Thus evaluating (\ref{bar-u2-tau2}) for $|x-q^k|=\tau$ and comparing it with (\ref{u2-tau}), we have

\begin{align}\label{2-q-1}
\bar u_2^k=-\lambda_{2,q}^k+2\log \delta_k-2\log \frac{h_2(q^k)}8-12\pi \beta(q^k,q^k)-8\pi \sum_{l\ge 3}G(q^k,p_l^k)\nonumber\\
+SE_3\log \delta_k^{-1}+SE_4.
\end{align}

Now we evaluate $u_2^k$ in a neighborhood of $p_l^k$ for each $l\ge 3$, from the Green's representation of $u_2^k$ we have
\begin{align}\label{upl-tem}
u_2^k(x)=&\bar u_2^k+12\pi G(x,q^k)-4\log |x-p_l^k|+8\pi \beta (x,p_l^k)+\sum_{s\neq l}8\pi G(x,p_s^k)\\
&+SE_4+SE_3 \log \delta_k^{-1}. \nonumber
\end{align}
Use the expansion of $u_2^k$ around $p_l^k$ and the result for single equation we have
\begin{align}\label{upl-tem-2}
\tilde u_2^k(x)&=u_2^k(x)-f_2^k(x)-\psi_2^k(x)\\
&=\log \frac{e^{\lambda_{2,l}^k}}{(1+e^{\lambda_{2,l}^k}\frac{h_2(p_l^k)}4|x-\tilde p_l^k|^2)^2}+SE_4 \nonumber
\end{align}
for some $\tilde p_l^k-p_l^k=O(e^{-\lambda_2^k})$. In the neighborhood of $p_l^k$ we have
\begin{align}\label{upl-tem-3}
&f_2^k(x)+\psi_2^k(x)
= 12\pi(G(x,q^k)-G(p_l^k,q^k))\\
&+8\pi(\beta(x,p_l^k)
-\beta(p_l^k,p_l^k))+8\pi\sum_{s\neq l}(G(x,p_s^k)-G(p_l^k,p_s^k)). \nonumber
\end{align}
Comparing (\ref{upl-tem}), (\ref{upl-tem-2}) and (\ref{upl-tem-3}) we have
\begin{align}
\bar u_2^k=-\lambda_{2,l}^k-2\log \frac{h_2(p_l^k)}4-12\pi G(p_l^k,q^k)-8\pi \beta (p_l^k,p_l^k)-\sum_{s\neq l,s\ge 3}8\pi G(p_l^k,p_s^k) \nonumber\\
+SE_4+SE_3\log \delta_k^{-1}.\label{l-u2k}
\end{align}
The two different expression of $\bar u_2^k$ in  (\ref{2-q-1}) and (\ref{l-u2k}) lead to
\begin{align}\label{lam2-q}
\lambda_{2,q}^k-2\log \delta_k+2\log \frac{h_2(q^k)}8+12\pi \beta(q^k,q^k)+8\pi \sum_{l\ge 3}G(q^k,p_l^k) \nonumber \\
=\lambda_{2,l}^k+2\log \frac{h_2(p_l^k)}4+12\pi G(p_l^k,q^k)+8\pi \beta (p_l^k,p_l^k)+\sum_{s\neq l,s\ge 3}8\pi G(p_l^k,p_s^k)\nonumber \\
+O(\delta_k)+O(\lambda_2^k\lambda_2^ke^{-\lambda_2^k}\delta_k^{-4-2d_k})+O(\lambda_1^ke^{-\lambda_1^k}\delta_k^{-4-d_k}\log \delta_k^{-1}).
\end{align}

Finally using results for single equation and (\ref{upl-tem}) we know that at each $p_l^k$ we have
\begin{align}\label{grad-pl-case-3}
\nabla (\log h_2)(p_l^k)+12\pi\nabla_1G(p_l^k,q^k)+8\pi\nabla_1\beta(p_l^k,p_l^k)+8\pi\sum_{s\neq l}\nabla_1 G(p_l^k,p_s^k)\nonumber\\
=SE_3.
\end{align}

\medskip

\noindent{\bf Proof of Theorem \ref{thm4}:} The largeness of $\lambda_1^k$ is proved in (\ref{mu-1k-big}), the fact that $\lambda_2^k+2\log \delta_k\to \infty$ is the nature of the selection process that determines the bubbling disks. The estimate of $\rho_1^k-4\pi$ is a combination of (\ref{tem-11}) and (\ref{tem-12}). The corresponding result on $\rho_2^k-8\pi$ is proved in (\ref{rho-2-q}) and the smallness of $\int_{M\setminus B(q^k,\tau)}\rho_2^kh_2e^{u_2^k}$ (which is $O(e^{-\lambda_2^k})$). In this case there is only one common blowup point $q^k$, whose location satisfies (\ref{loc-q-2}).  Theorem \ref{thm4} is established. $\Box$

\medskip

\noindent{\bf Proof of Theorem \ref{thm-case-3-N}:}. The estimate of $\rho_1^k-4\pi$ comes from (\ref{tem-11}) and (\ref{tem-13}). To estimate $\rho_2^k-4N\pi$, first we have
(\ref{rho-2-q}) that concerns the integation of $\rho_2^k h_2e^{u_2^k}$ near $q^k$. For $l\ge 3$, by result for single Liouville equation, we have
 $$\rho_{2,l}^k-4\pi=O(\lambda_2^ke^{-\lambda_2^k})$$  where $\rho_{2,l}^k=\int_{B(p_l^k,\tau)}\rho_2^k h_2e^{u_2^k}$ for small $\tau$. Then using $u_2^k=-\lambda_2^k+O(1)$ outside bubbling area, we obtain the desired estimate of $\rho_2^k-4N\pi$ in the statement. The location requirements of $q^k$ and $p_l^k$ are derived in (\ref{loc-q-f}) and (\ref{grad-pl-case-3}), the comparison of $\lambda_{2,q}^k$ with $\lambda_{2,l}^k$ is established in (\ref{lam2-q}). Theorem \ref{thm-case-3-N} is established. $\Box$

\section{Case two}
\begin{thm}\label{thm-case-2} If the case two of Theorem \ref{thm1} occurs, we have, for $N=2$, the following estimates:
\begin{equation}\label{lam-1-big}
\lambda_1^k-2\lambda_2^k\to \infty,
\end{equation}
\begin{equation}\label{rho-2}
\rho_2^k-8\pi = O(e^{-\lambda_2^k/2})+O(e^{-\lambda_1^k+2\lambda_2^k}).
\end{equation}
\begin{equation}\label{rho-1-case2}
\rho_1^k-4\pi=(D+o(1))e^{-\lambda_1^k+2\lambda_2^k}
\end{equation}
for some $D>0$.
Finally at $q^k$ we have
$$
\nabla(\log h_1+2\log h_2)(q^k)+24\pi \nabla_1\beta(q^k,q^k)
=O(e^{-\lambda_2^k/2})+O(e^{-\lambda_1^k+2\lambda_2^k}).
$$
\end{thm}
Note that in case two, $q^k$ is the common blowup point for both $u_1^k$ and $u_2^k$. $u_2^k$ has no other blowup points if $N=2$.
\begin{thm}\label{thm-case-2N} In case two, if $N>2$,
\begin{equation}\label{lam-1-big-2}
\lambda_1^k-\lambda_2^k\to \infty, \quad \lambda_{2,q}^k=\frac 12 \lambda_2^k+O(1),
\end{equation}
\begin{equation}\label{rho-2N}
\rho_2^k-4\pi N
=O(e^{-\lambda_2^k/4})+O(e^{-\lambda_1^k+\lambda_2^k}).
\end{equation}
\begin{equation}\label{rho-1-case2N}
\rho_1^k-4\pi=(D+o(1))e^{-\lambda_1^k+\lambda_2^k}\int_M h_1 exp(-4\pi \sum_{l=3}^N G(x,p_l^k))dV_g.
\end{equation}
for some $D>0$.
Finally the locations of blowup points are related by
\begin{align*}
&\nabla (\log h_2)(p_l^k)+12\pi \nabla_1 G(p_l^k,q^k)+8\pi\sum_{s\neq l,s=3}^N\nabla_1 G(p_l^k,p_s^k)+8\pi \nabla_1\beta(p_l^k,p_l^k)\\
=&O(e^{-\lambda_2^k/4})+O(e^{-\lambda_1^k+\lambda_2^k}),\quad l=3,..,N\\
&\nabla (\log h_1+2\log h_2)(q^k)+24\pi \nabla_1\beta(q^k,q^k)+12\pi \sum_{l\ge 3}\nabla_1 G(q^k,p_l^k)\nonumber \\
=&O(e^{-\lambda_2^k/4})+O(e^{-\lambda_1^k+\lambda_2^k}).
\end{align*} and the magnitudes of $u_2^k$ at $p_l^k$ ($l\ge 3$) are linked by
\begin{align*}
\lambda_{2,l}^k+2\log \frac{h_2(p_l^k)}8+12 \pi G(p_l^k,q^k)+8\pi \sum_{s=3,s\neq l}^N G(p_l^k,p_s^k)
+8\pi \beta(p_l^k,p_l^k)\nonumber\\
=\lambda_{2,m}^k+2\log \frac{h_2(p_m^k)}8+12 \pi G(p_m^k,q^k)+8\pi \sum_{s=3,s\neq m}^N G(p_m^k,p_s^k)
+8\pi \beta(p_m^k,p_m^k) \nonumber \\
+O(e^{-\lambda_2^k/4})+O(e^{-\lambda_1^k+\lambda_2^k}),\quad \mbox{ for }\quad l,m=3,...,N \quad \mbox{ and } \quad l\neq m.
\end{align*}
\end{thm}
Theorem \ref{thm-case-2} immediately implies that case two does not occur if $\rho_1^k$ tends to $4\pi$ from below.
In local coordinates of $q^k$ we rewrite the equation just like in case three: let $\phi$ satisfy (\ref{phi-cur}), $f_1^k$ and $f_2^k$ be defined as in (\ref{case3-f}) and let $\psi_i^k$, $\tilde u_i^k$, $h_i^k$ be defined as in case three. Then in $B(0,\tau)$ we have (\ref{sys-q}) in $B(0,\tau)$.

The difference between this case and case three is that if we set $\epsilon_k=e^{-\lambda_{2,q}^k/2}$ and
$$v_2^k(y)=\tilde u_2^k(\epsilon_ky)+2\log \epsilon_k.  $$

$v_2^k$ converges to a function $v(y)$ that satisfies
$$\Delta v(y)+2\lim_{k\to \infty}  h_2^k(0)e^{v}=4\pi \delta_0, $$
where $h_1^k(\epsilon_k y)e^{v_1^k(y)}$ tends to $4\pi \delta_0$ in measure. Here
$$v_1^k(y)=\tilde u_1^k(\epsilon_ky)+2\log \epsilon_k.$$
So in this case we have
$$\lambda_1^k-\lambda_{2,q}^k\to \infty.$$
Let $\Omega_k=B(0,\tau\epsilon_k^{-1})$
and we use the following function to remove $v_2^k$ from the equation for $v_1^k$: $f_{2b}^k$ defined by
$$f_{2b}^k(y)=\frac 1{2\pi}\int_{\Omega_k}\log |y-\eta |  h_2^k(\epsilon_k \eta)e^{v_2^k(\eta)}d\eta, $$
satisfies
$$\Delta f_{2b}^k(y)= h_2^k(\epsilon_k y)e^{v_2^k(y)}.$$
Thus we write the equation for $v_1^k$ as
\begin{equation}\label{tem-case-2}
\Delta (v_1^k-f_{2b}^k)+2 h_1^k(\epsilon_k y)e^{f_{2b}^k}e^{v_1^k-f_{2b}^k}=0.
\end{equation}
Here
we further observe that
\begin{equation}\label{mu-1-large}
\mu_1^k+2\log \epsilon_k\to \infty,\quad \mu_1^k:=\lambda_1^k+2\log \epsilon_k.
\end{equation}
The reason is for $v_1^k$, $v_2^k$ in $\Omega_k\setminus B_R$ for $R$ large. The fast decay of $v_2^k$ and the fact that $v_2^k\le 0$ makes it easy to use a test function to remove $v_2^k$ in the equation for $v_1^k$.  Then just like the proof in the previous section, starting from $\bar v_1^k(R)=-\mu_1^k+O(R)$, we can use $\rho_1^k\to 4\pi$ to further prove (\ref{mu-1-large}) just like before.

Clearly (\ref{mu-1-large}) is equivalent to
\begin{equation}\label{lam-1-large}
\lambda_1^k-2\lambda_{2,q}^k\to \infty.
\end{equation}

The derivation of the leading term of $\rho_1^k-4\pi$ is similar to case three. After scaling according to $\lambda_{2,q}^k$, the total integration of $\rho_1^kh_1e^{u_1^k}$ can be estimated by two parts. One part is the integration in $B(q^k,R\epsilon_k)$ where $\epsilon_k=e^{-\lambda_{2,q}^k/2}$ , $R>>1$ and the remaining part is the integration over $M\setminus B(q^k, R\epsilon_k)$. Using Gluck's result in \cite{gluck} for (\ref{tem-case-2}) and the fact that $\Delta f_{2b}^k(0)=O(1)$, we have
\begin{equation}\label{rho-1-case-2}
\int_{B(q^k,R\epsilon_k)}\rho_1^kh_1e^{u_1^k}dV_g
=4\pi+O(\mu_1^ke^{-\mu_1^k}\epsilon_k^2)+O(e^{-\mu_1^k}).
\end{equation}
  If $N=2$, by the Green's representation of $u_1^k$, we see that for $x\not \in B(q^k,\tau)$,
\begin{align*}
u_1^k(x)&=\bar u_1^k+O(\mu_1^ke^{-\lambda_1^k})+O(e^{-\mu_1^k+\lambda_2^k})\\
&=\bar u_1^k+O(e^{-\mu_1^k+\lambda_2^k}), \mbox{ by (\ref{both-b})},
\end{align*}
which implies that
$$\bar u_1^k=-\mu_1^k-2\log \epsilon_k+O(1)=-\lambda_1^k-4\log \epsilon_k+O(1).$$
 Thus
\begin{equation}\label{rho-1-case-2b}
\int_{M\setminus B(q^k,R\epsilon_k)}\rho_1^kh_1 e^{u_1^k}dV_g=(D+o(1))e^{-\lambda_1^k+2\lambda_{2,q}^k}
\end{equation}
and in this case $\lambda_{2,q}^k=\lambda_2^k$.

If $N\ge 3$ and $x\not \in B(q^k,\tau)$,
$$u_1^k(x)=\bar u_1^k-4\pi \sum_{l\ge 3} G(x,p_l^k)+O(\lambda_2^ke^{-\lambda_2^k})+O(e^{-\lambda_1^k+\lambda_{2,q}^k}). $$
Thus in this case
\begin{equation}\label{rho-1-case-2c}
\int_{M\setminus B(q^k,R \epsilon_k)} \rho_1^k h_1e^{u_1^k}dV_g=(D+o(1))\int_M h_1 e^{-\sum_l 4\pi G(x,p_l^k)}dV_g e^{-\lambda_1^k+\lambda_2^k}.
\end{equation}

The combination of (\ref{rho-1-case-2}) and (\ref{rho-1-case-2b}) leads to
(\ref{rho-1-case2}) in Theorem \ref{thm-case-2}. (\ref{rho-1-case2N}) in Theorem \ref{thm-case-2N} is based on (\ref{rho-1-case-2}) and (\ref{rho-1-case-2c}).

Here we further remark that the difference between $\lambda_1^k$ and $\lambda_{2,q}^k$ in (\ref{lam-1-large}) implies
\begin{equation}\label{case-2-v2}
v_2^k(y)=-6\log |y|+O(1),\quad |y|\sim \epsilon_k^{-1}.
\end{equation}
The proof of (\ref{case-2-v2}) is very similar to the corresponding estimate for case three. The fast decay of $v_2^k$ and the largeness of $\mu_1^k$ ((\ref{mu-1-large}) yield the accurate estimate of $v_2^k$ in (\ref{case-2-v2}).

Before we compute the leading term of the $\rho_2^k-4\pi N$ we first observe that
\begin{equation}\label{q-2-1}
\lambda_{2,q}^k=\frac 12 \lambda_{2,l}^k+O(1),\quad l=3,...,N
\end{equation}

Indeed, according to (\ref{case-2-v2}), $u_2^k=-2\lambda_{2,q}^k+O(1)$ on $\partial B(q^k,\tau)$. On the other hand $u_2^k(x)=-\lambda_{2,l}^k+O(1)$ on
$\partial B(p_l^k,\tau)$ for $\tau>0$ small. Thus (\ref{q-2-1}) holds.

Before we compute $\rho_2^k-4\pi N$, we treat $h_1^k(\epsilon_ky)e^{v_1^k}$ as a perturbation of $4\pi \delta_0$ in the equation for $v_2^k$.

Let $f_{1b}^k(y)$ be defined by
$$f_{1b}^k(y)=\frac 1{2\pi }\int_{\Omega_k} \log |y-\eta |h_1^k(\epsilon_k \eta)e^{v_1^k(\eta )}d\eta, \quad y\in \Omega_k $$
where $\Omega_k=B(0,\tau \epsilon_k^{-1})$.
Clearly $f_{1b}^k$ solves
$$\Delta f_{1b}^k(y)=h_1^k(\epsilon_ky)e^{v_1^k}\quad \mbox{in}\quad \Omega_k. $$
The following lemma proves a minor difference between $e^{f_{1b}^k}$ and $|y|^2$.

\begin{lem}\label{extra-s}
\begin{equation}\label{one-est}
e^{f_{1b}^k(y)}-|y|^2=
O(e^{-\mu_1^k+\lambda_{2,q}^k})\log (2+|y|),\quad y\in \Omega_k.
\end{equation}
\end{lem}

\noindent{\bf Proof of Lemma \ref{extra-s}:}
We evaluate $f_{1b}^k$ in two parts:
\begin{equation}\label{f1b-t3}
f_{1b}^k(y)=\frac{1}{2\pi}\log |y| \int_{\Omega_k}h_1^k(\epsilon_k\eta)e^{v_1^k(\eta)}d\eta+\frac{1}{2\pi}\int_{\Omega_k}\log \frac{|y-\eta|}{|y|}h_1^k(\epsilon_k\eta)e^{v_1^k(\eta)}d\eta.
\end{equation}
The first term is
$$(2+O(\mu_1^ke^{-\lambda_1^k})+O(e^{-\mu_1^k})+O(e^{-\mu_1^k+\lambda_{2,q}^k}))\log |y| $$
as a result of integration over $B(q^k,\tau)$ and $M\setminus B(q^k,\tau)$.
The evaluation of the second term is based on integration over the following three regions:
\begin{align*}
&\Sigma_1:=\{\eta\in \Omega_k;\quad |\eta |<|y|/2\},\\
&\Sigma_2:=\{\eta\in \Omega_k;\quad |\eta -y|<|y|/2\},\\
& \Sigma_3=\Omega_k\setminus (\Sigma_1\cup\Sigma_2).
\end{align*}

Apply standard estimate in each of the domains we see the second term is $O(e^{-\mu_1^k})(1+|y|)^2\log (2+|y|)$.
 Lemma \ref{one-est} follows from the combinations of estimates for these two parts. $\Box$

\medskip

\begin{rem}\label{f1b-bry}
For later use we observe from (\ref{f1b-t3}) and the estimate of $f_{1b}^k$ that the oscillation of $f_{1b}^k$ on $\partial \Omega_k$ is $O(e^{-\mu_1^k+\lambda_{2,q}^k})\epsilon_k^2\log \epsilon_k^{-1}$.
\end{rem}

Let $v_k=v_2^k-f_{1b}^k$, we write the equation for $v_k$ as
$$\Delta v_k+2|y|^2 h_2^k(\epsilon_ky)e^{v_k}=E_k $$
where
$$E_k=(|y|^2-e^{f_{1b}^k})h_2^k(\epsilon_ky)e^{v_2^k}=(O(\epsilon_k^2)+O(e^{-\mu_1^k}\epsilon_k^{-2}))(1+|y|)^{-6}\log(2+|y|). $$
$v_k$ converges to $U$ that satisfies
$$\Delta U+2\lim_{k\to \infty} h_2^k(0)|y|^2e^U=0,\quad \mbox{in}\quad \mathbb R^2 $$
with the expression
$$U(y)=\log \frac{\Lambda}{2(1+\frac{\Lambda}{32}|y^2-\xi|^2)^2 \lim_{k\to \infty} h_2^k(0)} $$
for some $\xi\in \mathbb R^2$ and $\Lambda>0$.

Using the idea of \cite{lwz-jems}, we can adjust the parameters of $U$ a little bit to $U_k$ and make $v_k$ agree with $U_k$ on three points $p_1$,$p_2$,$p_3$:
$$U_k=\log \frac{\Lambda_k}{2(1+\frac{\Lambda_k}{32}|y^2-\xi_k|^2)^2  h_2^k(0)} $$
where $\Lambda_k\to \Lambda$, $\xi_k\to \xi$.
Here we remark that $U_k$ is a solution of
$$\Delta U_k+2h_2^k(0)|y|^2e^{U_k}=0,\quad \mbox{in}\quad \mathbb R^2, \quad 2h_2^k(0)\int_{\mathbb R^2}|y|^2e^{U_k}=16\pi. $$

The detail is as follows:
Choose $1<<|p_1|<<|p_2|<<|p_3|$ such that the following matrix invertible (see the appendix for the details of choosing $p_i$):
\begin{equation}\label{ap-ma}
\left(\begin{array}{ccc}
\frac{\partial U}{\partial \Lambda}(p_1) & \frac{\partial U}{\partial \Lambda}(p_2) & \frac{\partial U}{\partial \Lambda}(p_3)\\
\frac{\partial U}{\partial \xi_1}(p_1) & \frac{\partial U}{\partial \xi_1}(p_2) & \frac{\partial U}{\partial \xi_1}(p_3)\\
\frac{\partial U}{\partial \xi_2}(p_1) & \frac{\partial U}{\partial \xi_2}(p_2) & \frac{\partial U}{\partial \xi_2}(p_3)
\end{array}
\right)
\end{equation}
where $\xi=\xi_1+i\xi_2$. Thus if a $o(1)$ perturbation is placed on $v$ (to make $v_k(p_j)=U_k(p_j)$ for $j=1,2,3$), all we need to do is change the parameters $\Lambda$, $\xi$ by a comparable amount. So even though we have a sequence of parameters $\Lambda_k$, $\mu_k$, they are not tending to infinity.
Now we rewrite the equation for $v_k$ as
\begin{equation}\label{vk-g}
\Delta v_k+2|y|^2h_2^k(0)e^{v_k}
=-2\epsilon_k\sum_{t=1}^2\partial_t h_2^k(0)y^t|y|^2e^{v_k}+\tilde E_k
\end{equation}
where the error term $\tilde E_k$ satisfies
$$|\tilde E_k|\le C(\epsilon_k^2+O(e^{-\mu_1^k+\lambda_{2,q}^k}))(1+|y|)^{-4}. $$
Let
$$w_k(y)=v_k(y)-U_k(y),\quad y\in \Omega_k:=B(0,\tau \epsilon_k^{-1}). $$
We have
$$w_k(y)=O(1) \quad \mbox{in}\quad \Omega_k $$
and
$$w_k(y)=C+O(\epsilon_k^2)+O(e^{-\mu_1^k+\lambda_{2,q}^k}) \log \epsilon_k^{-1}\quad \mbox{on}\quad \partial \Omega_k $$
because $v_2^k$ is a constant on $\partial \Omega_k$, $f_{1b}^k$ has an oscillation of $O(e^{-\lambda_1^k+\lambda_{2,q}^k}\lambda_{2,q}^k)$ ( see Remark \ref{f1b-bry}) and by the following expression of $U_k$,
the oscillation of $U_k$ on $\partial \Omega_k$ is $O(\epsilon_k^2)$:
\begin{align}\label{asy-U}
& U_k(y)
=-8\log r-\log (2h_2^k(0))-\log \Lambda_k+2\log 32\\
&+\frac 1{r^2}e^{2i\theta}\bar \xi_k\frac{32}{\Lambda_k}
+\frac 1{r^2}\xi_k e^{-2i\theta}\frac{32}{\Lambda_k}+O(r^{-4})\nonumber \\
&=-8\log r-\log \Lambda_k+2\log 32-\log (2h_2^k(0))\nonumber\\
&+\frac{64}{\Lambda_k}r^{-2}(\cos 2\theta \xi_1^k+\sin 2\theta \xi_2^k)+O(r^{-4}). \nonumber
\end{align}

Moreover, we have
$$w_k(p_j)=0,\quad j=1,2,3. $$
The equation for $w_k$ is
$$\Delta w_k+2|y|^2h_2^k(0)e^{\xi_k}w_k=-2\epsilon_k\sum_t\partial_th_2^k(0)y^t|y|^2e^{v_k}+\tilde E_k $$
where $\xi_k$ comes from the Mean Value Theorem.  The main result for $w_k$ is
\begin{lem}\label{case-2-est-1}Let $\bar \epsilon_k=\epsilon_k+e^{-\mu_1^k+\lambda_{2,q}^k}$, then
\begin{equation}\label{crude-wk}
|w_k(y)|\le C\bar \epsilon_k\log(2+|y|),\quad y\in \Omega_k.
\end{equation}
\end{lem}

\noindent{\bf Proof of Lemma \ref{case-2-est-1}}:
The proof of (\ref{crude-wk}) is by contradiction. Suppose
$$\tilde \Lambda_k:=\max_{y\in \bar\Omega_k}\frac{|w_k(y)|}{(\epsilon_k+e^{-\mu_1^k+\lambda_{2,q}^k})\log (2+|y|)}\to \infty$$
and $\tilde \Lambda_k$ is attained at $y_k$ (which could appear on $\partial \Omega_k$). We set
$$\hat w_k(y)=\frac{w_k(y)}{\tilde \Lambda_k(\epsilon_k+e^{-\mu_1^k+\lambda_{2,q}^k})\log (2+|y_k|)}, $$
which is obviously a solution of
$$\Delta \hat w_k(y)+2|y|^2h_2^k(0)e^{\xi_k}\hat w_k=\frac{O(1+|y|)^{-4}}{\tilde \Lambda_k \log (2+|y_k|)}$$
and the definition of $\hat w_k$ implies
$$|\hat w_k(y)|\le \frac{\log (2+|y|)}{\log (2+|y_k|)}. $$

Green's representation formula for $\hat w_k(y_k)$ gives
\begin{align*}
&\pm 1=\hat w_k(y_k)=\hat w_k(y)-\hat w_k(p_1)\\
&=\int_{\Omega_k}(G(y_k,\eta)-G(p_1,\eta))(2h_1^k(\epsilon_k\eta)|\eta |^2 e^{\xi_k}\hat w_k(\eta)+\frac{O(1)(1+|\eta|)^{-4}}{\tilde \Lambda_k\log (2+|y_k|)})d\eta\\
&+\frac{O(1)}{\tilde \Lambda_k\log (2+|y_k|)}.
\end{align*}
where we have used $\hat w_k(p_1)=0$ and $\hat w_k=C+\frac{O(\epsilon_k)}{\tilde \Lambda_k\log (2+|y_k|)}$ on $\partial \Omega_k$.

The estimate of the Green's function $G_k$ on $\Omega_k$ is (see \cite{lin-zhang-jfa} for detail)
$$|G(y_k,\eta)-G(p_1,\eta)|\le \left\{\begin{array}{ll}
C(\log |\eta |+\log |y_k|),\quad \eta\in \Sigma_1,\\
C(\log |y|+|\log |y-\eta||,\quad \eta \in \Sigma_2,\\
C|y|/|\eta |\quad \eta\in \Sigma_3.
\end{array}
\right.
$$
where
\begin{align*}
&\Sigma_1=\{\eta\in \Omega_k;\quad |\eta|<|y|/2\quad \}\\
&\Sigma_2=\{\eta\in \Omega_k;\quad |\eta -y|<|y|/2,\quad \}\\
&\Sigma_3=\Omega_k\setminus (\Sigma_1\cup \Sigma_2).
\end{align*}
If $y_k\to y^*$, $\hat w_k$ converges to a solution of
$$\Delta \phi+2\lim_{k\to \infty}h_2^k(0)|y|^{2}e^{U}\phi=0,\quad \mathbb R^2, $$
with mild growth:
$$|\phi(y)|\le C\log (2+|y|). $$
By the non-degeneracy of the linearized equation,
$$\phi(y)=c_1 \frac{\partial U}{\partial \Lambda}(y)+c_2\frac{\partial U}{\partial \xi_1}(y)+c_3\frac{\partial U}{\partial \xi_2}(y).$$
Using $\phi(p_i)=0$ for $i=1,2,3$, we have, by the invertibility of matrix (\ref{ap-ma}), $c_1=c_2=c_3=0$, thus $\phi\equiv 0$, a contradiction to
$\hat w_k(y_k)=\pm 1$.

If $y_k\to \infty$, the evaluation of $\hat w_k(y_k)=o(1)$ can be obtained by elementary estimates, a contradiction to $\pm 1=\hat w_k(y_k)$. Thus (\ref{crude-wk}) is established. $\Box$

\medskip

Just as in case three, the estimate of $\rho_2^k-4N\pi$ starts from the neighborhood of $q^k$: let $\tilde \rho_2^k=\rho_2^k\int_{B(q^k,\tau)}h_2e^{u_2^k}dV_g$ and
$\bar \epsilon_k=\epsilon_k+e^{-\lambda_1^k+2\lambda_{2,q}^k}$, then
\begin{align}\label{case-2-rho-22}
2\tilde \rho_2^k&=\int_{B(q^k,\tau)}2\rho_2^kh_2e^{u_2^k}dV_g\\
&=2\int_{\Omega_k}|y|^2 h_2^k(\epsilon_ky)e^{v_2^k(y)}dy+O(\epsilon_k^2) \quad \mbox{where }\quad \Omega_k=B(0,\tau \epsilon_k^{-1})\nonumber\\
&=2\int_{\Omega_k}|y|^2 ( h_2^k(0)+O(\epsilon_k)|y|))e^{U_k}(1+O(\bar\epsilon_k)\log (2+|y|)dy\nonumber\\
&=16\pi+O(\bar \epsilon_k).\nonumber
\end{align}

From here we further claim that, for $N\ge 3$,
\begin{equation}\label{rho-2-case-2N}
\rho_2^k-4\pi N=O(\bar \epsilon_k).
\end{equation}
Indeed, besides the energy in $B(q^k,\tau)$, we observe that the total integration of $\rho_2^kh_2e^{u_2^k}$ away from bubbling disks is $O(e^{-\lambda_2^k})$. The integration around each $p_l$ for $l\ge 3$ is $O(\lambda_2^ke^{-\lambda_2^k})$ by results for single equation. Thus (\ref{rho-2-case-2N}) holds and (\ref{rho-2N}) in Theorem \ref{thm-case-2N} is verified.

\medskip

Now we discuss the location of blowup points and the comparison of bubble heights in case two.
For convenience of notation we set
$$SE_5=O(e^{-\lambda_{2,q}^k/2})+O(\lambda_1^k e^{-\lambda_1^k+\lambda_{2,q}^k})+O(e^{-\lambda_1^k+2\lambda_{2,q}^k})=O(e^{-\lambda_{2,q}^k/2})+O(e^{-\lambda_1^k+2\lambda_{2,q}^k}), $$
which is the main error in this section.
First we use  the Pohozaev identity for system to determine the location of $q^k$:  For $\xi\in \mathbb S^1$,
\begin{align}\label{pi-sys}
&\sum_i\int_{B_{\tau}}\partial_{\xi}(\log h_i^k)h_i^ke^{\tilde u_i^k}\\
=&\sum_i\int_{\partial B_{\tau}}h_i^ke^{\tilde u_i^k}(\xi\cdot \nu)+a^{ij}\partial_{\nu}\tilde u_j^k\partial_{\xi}\tilde u_i^k-\frac 12
a^{ij}(\nabla \tilde u_i^k\nabla \tilde u_j^k)(\xi\cdot \nu)dS, \nonumber
\end{align}
where $\nu$ is the outer-normal vector on $\partial B_{\tau}$.

By using the expansion of $u_1^k$ and $u_2^k$ at $q^k$, we see that
\begin{equation}\label{q-case-2}
\nabla \log h_1^k(0)+2\nabla \log h_2^k(0)=SE_5.
\end{equation}

Note that (\ref{q-case-2}) is obtained as follows. First for $\int_{B_{\tau}}\partial_{\xi}(\log h_1^k)h_1^ke^{\tilde u_1^k}$ we use the point-wise analysis for $\tilde u_1^k$ inside $B(0, N\epsilon_k)$ and $B_{\tau}\setminus B(0, N \epsilon_k)$ for $N$ large. Then it is easy to see that
\begin{align*}
\int_{B_{\tau}}\partial_{\xi}(\log h_1^k)h_1^ke^{\tilde u_1^k}&=4\pi \partial_{\xi}(\log h_1^k)(0)+O(\mu_1^ke^{-\lambda_1^k})+O(e^{-\lambda_1^k+2\lambda_{2,q}^k})\\
&=4\pi \partial_{\xi}(\log h_1^k)(0)+O(e^{-\lambda_1^k+2\lambda_{2,q}^k}).
\end{align*}

The integration of $\partial_{\xi}(\log h_2^k)h_2^ke^{\tilde u_2^k}$ can be evaluated using the $O(\bar \epsilon_k)$ expansion of $\tilde u_2^k$:
$$\int_{B_{\tau}}\partial_{\xi}(\log h_2^k)h_2^ke^{\tilde u_2^k}=8\pi \partial_{\xi}(\log h_2^k)(0)+O(\bar \epsilon_k). $$
Since $\xi$ is an arbitrary vector in $\mathbb S^1$, the left hand side of (\ref{pi-sys}) is
$$4\pi\partial_{\xi}\log h_1^k(0)+8\pi \partial_{\xi}h_2^k(0)+SE_5. $$
It is also easy to evaluate the right hand side of (\ref{pi-sys}), because, first, $\int_{\partial B_{\tau}}h_i^ke^{u_i^k}(\xi\cdot \nu)=SE_5$ by the rough estimate of $u_i^k$ on $\partial B(q^k,\tau)$. For the two other terms on the right hand side, we use Green's representation of $\tilde u_i^k$ on $B_{\tau}$. Since the value of $\tilde u_i^k$ is a constant on $\partial B_{\tau}$, it disappears after differentiation. Then
\begin{align*}
\nabla \tilde u_1^k(x)=SE_5\\
\nabla \tilde u_2^k(x)=-6\frac{x}{|x|^2}+SE_5.
\end{align*}
Using these expressions in the evaluation of the two remaining terms, the right hand side of the Pohzoaev identity is $SE_5$.  Thus (\ref{q-case-2}) holds.

To further determine (\ref{q-case-2}) we consider $N=2$ and $N\ge 3$ separately.
First for $N=2$, the Green's representation of $u_1^k$ gives
$$u_1^k(x)=\bar u_1^k+\int_M G(x, \eta)(2\rho_1^kh_1^ke^{u_1^k}-\rho_2^kh_2e^{u_2^k})dV_g(\eta), $$
and by evaluating $u_1^k$ away from $q^k$ it is easy to obtain
$$u_1^k(x)=\bar u_1^k+SE_5, \quad x\in M\setminus B(q^k,\tau), $$
which gives
\begin{align*}
\nabla_1\psi_1^k(q^k)=SE_5,\quad \mbox{if}\quad N=2, \\
\nabla_1\psi_1^k(q^k)=12\pi \nabla_1\beta(q^k,q^k)+SE_5,\quad \mbox{if}\quad N\ge 3.
\end{align*}
Thus (\ref{q-case-2}) is translated as the following for $N=2$:
\begin{equation}\label{q-case-2-2}
\nabla(\log h_1+2\log h_2)(q^k)+24\pi \nabla_1\beta(q^k,q^k)=SE_5.
\end{equation}

\medskip

Next for $N\ge 3$ we not only determine the location of $q^k$, but also $p_l^k$ for $l\ge 3$.
Green's representation of $u_2^k$ gives
$$u_2^k(x)=\bar u_2^k+\int_M G(x,\eta)(2\rho_2^kh_2e^{u_2^k}-\rho_1^kh_1 e^{u_1^k})dV_g(\eta), $$
where $\bar u_2^k$ is the average of $u_2^k$.
After using results for single equation we have
$$u_2^k(x)=\bar u_2^k+\sum_{s=3}^N8\pi G(x, p_s^k)+12\pi G(x,q^k)+SE_5. $$

By the same method for case one, we see that the locations must satisfy
\begin{align}\label{case-2-loc}
\nabla (\log h_2)(p_l^k)+12\pi \nabla_1 G(p_l^k,q^k)+8\pi\sum_{s\neq l,s=3}^N\nabla_1 G(p_l^k,p_s^k)+8\pi \nabla_1\beta(p_l^k,p_l^k) \nonumber\\
=SE_5.
\end{align}

For $N\ge 3$ we have a different expression of the expansions of $u_i^k$ around $q^k$:
$$u_1^k(x)=\bar u_1^k-\sum_{l\ge 3}4\pi G(x,p_l^k)+O(\lambda_2^ke^{-\lambda_2^k}), $$
$$u_2^k(x)=\bar u_2^k+\sum_{l\ge 3}8\pi G(x,p_l^k)+12\pi G(x,q^k)+SE_5. $$
So the derivatives of harmonic functions that cancel the oscillation of $u_i^k$ around $q^k$ are
\begin{align}\label{q-case-2-N1}
\nabla\psi_1^k(q^k)=-4\pi \sum_{l\ge 3}\nabla_1G(q^k,p_l^k)+SE_5 \\
\nabla\psi_2^k(q^k)=12\pi \nabla_1\beta(q^k,q^k)+\sum_{l\ge 3}8\pi \nabla_1G(q^k,p_l^k)+SE_5. \nonumber
\end{align}
Using (\ref{q-case-2-N1}) in (\ref{q-case-2}), we have
\begin{align}\label{loc-q-case-2N}
\nabla (\log h_1+2\log h_2)(q^k)+12\pi \sum_{l\ge 3}\nabla_1 G(q^k,p_l^k)+24\pi \nabla_1\beta(q^k,q^k)\nonumber \\
=SE_5.
\end{align}

Finally we compare the magnitudes of $u_2^k$ at $p_l^k$  for $l\ge 3$.
By evaluating $u_2^k(p_l^k)$ using standard estimates for single Liouville equation, we have
\begin{align*}
-\bar u_2^k=\lambda_{2,l}^k+\frac 2{h_2(p_l^k)}\log \frac{h_2^k(p_l^k)}8+12 \pi G(p_l^k,q^k)+8\pi \sum_{s=3,s\neq l}^N G(p_l^k,p_s^k)\\
+8\pi \beta(p_l^k,p_s^k)+O(\lambda_2^k\lambda_2^k e^{-\lambda_2^k/2}).
\end{align*}

Equating $\bar u_2^k$ for $m\neq l$, we arrive at

\begin{eqnarray}
\lambda_{2,l}^k+2\log \frac{h_2(p_l^k)}8+12 \pi G(p_l^k,q^k)+8\pi \sum_{s=3,s\neq l}^N G(p_l^k,p_s^k)
+8\pi \beta(p_l^k,p_l^k)\nonumber\\
\lambda_{2,m}^k+2\log \frac{h_2(p_m^k)}8+12 \pi G(p_m^k,q^k)+8\pi \sum_{s=3,s\neq m}^N G(p_m^k,p_s^k)
+8\pi \beta(p_m^k,p_m^k) \nonumber \\
+O(\lambda_2^k\lambda_2^ke^{-\lambda_2^k/2})+SE_5. \label{comp-2}
\end{eqnarray}

These are about the location of $p_l^k$ for $l\ge 3$.

\medskip

\noindent{\bf Proof of Theorem \ref{thm-case-2}:} The comparison of $\lambda_1^k$ and $\lambda_2^k$, which is $\lambda_{2,q}^k$ in this case is stated in (\ref{lam-1-large}). The combination of (\ref{rho-1-case-2}) and (\ref{rho-1-case-2b}) leads to
(\ref{rho-1-case2}). The difference between $\rho_2^k$ and $8\pi$ is implied by (\ref{case-2-rho-22}). Finally the location of $q^k$ is established in (\ref{q-case-2-2}). Theorem \ref{thm-case-2} is established. $\Box$

\medskip

\noindent{\bf Proof of Theorem \ref{thm-case-2N}}: The relation between $\lambda_1^k$ and $\lambda_2^k$, as well as $\lambda_{2,q}^k$ is stated in (\ref{lam-1-large}) and (\ref{q-2-1}).
The estimate of $\rho_1^k-4\pi$ is a combination of (\ref{rho-1-case-2}) and (\ref{rho-1-case-2c}). The difference between $\rho_2^k-4N$, stated in (\ref{rho-2N}), is verified by (\ref{rho-2-case-2N}). The locations of $q^k$ and $p_l^k$ can be found in (\ref{loc-q-case-2N}) and (\ref{case-2-loc}). Finally the comparison of local maximums of $u_2^k$ other than $q^k$ is stated in (\ref{comp-2}).
Theorem \ref{thm-case-2N} is established. $\Box$

\section{Proof of theorems stated in the introduction}

\noindent{\bf Proof of Theorem \ref{thm-simple}:}  This is immediate for $\rho_1^k\to 4\pi$ and $\rho_2^k\to 4\pi$. Since the classification of $(\sigma_1^k,\sigma_2^k)$ does not allow the bubbling disk of $u_1^k$ to collide with that of $u_2^k$.
Theorem \ref{thm-simple} is established $\Box$

\medskip

\noindent{\bf Proof of Theorem \ref{thm1}:}  From the classification of $(\sigma_1(q^k),\sigma_2(q^k))$ we see that $(2,4)$ is the only case that bubbling disk collision is possible. Case three and case two have been discussed. Here we further claim that there is no other case. Here we put into two cases, if the spherical Harnack inequality around $q^k$ holds for $u_2^k$, we prove that case two happens. If the spherical Harnack does not hold we prove that case three happens. First if the spherical Harnack holds around $q^k$ for $u_2^k$, it is not possible to have $\lambda_{2,q}^k-\lambda_1^k\to \infty$, because otherwise the scaling of $u_1^k$ according to its maximum will have singular sources, which will make $\rho_1^k$ tend to a value greater than $4\pi$. If $\lambda_1^k=\lambda_{2,q}^k+O(1)$, one sees immediately that this is not possible. First the bubbling profile in a disk around $q^k$ already has $\lambda_1^k+2\log \delta_k\to \infty$, where $\delta_k\to 0$ in the radius of bubbling disk around $q^k$, in which a profile of global solution of $u_1^k$ can been seen inside. If $\lambda_{2,q}^k\ge \lambda_1^k+O(1)$ we would have
$\lambda_{2,q}^k+2\log \delta_k\to \infty$, which is case three instead of case two. Therefore for case two we must have $\lambda_1^k-\lambda_{2,q}^k\to \infty$. The discussion for case two further gives
$\lambda_1^k-2\lambda_{2,q}^k\to \infty$.

For case three, we here remark that it is not possible to have two bubbling disks of $u_2^k$ not in comparable distance to $q^k$: In other words, it is not possible to have $p_1^k,p_2^k$ both tending to $q$ such that $|p_1^k-q^k|/|p_2^k-q^k|\to \infty$. The reason is if this happen, using the argument in \cite{lwyz-apde} and \cite{lwz-apde}, there is $l_k\to 0$ such that Pohozaev identity can be evaluated on $B(q^k,l_k)$ but $l_k/|p_1^k-q^k|\to 0$, $l_k/|p_2^k-q^k|\to \infty$
and $u_2^k$ and $u_1^k$ both have fast decay on $B(q^k,l_k)$, which means Pohozaev identity gives $(2,2)$ as a type of concentration of $(\sigma_1(q),\sigma_2(q))$, which is against the known classification result for $(\sigma_1(q),\sigma_2(q))$. Then the further discussions about case two and case three finished the proof of Theorem \ref{thm1}. $\Box$

\medskip

\noindent{\bf Proof of Theorem \ref{simple-2}:} First from (\ref{rho-1-case2}) in Theorem \ref{thm-case-2}  we see that case two does not occur. In order to rule out case three, we  observe from the
 expansion of $\rho_1^k-4\pi$ that the only way that $\rho_1^k<4\pi$ is when $\delta_k^{-4+\epsilon}$ is less than $\lambda_1^k$. Then, in view of (\ref{both-b}), it is enough to see $e^{-\lambda_2^k}=o(\delta_k^{-12})$. If this happens, we see immediately that
$$\nabla (\log h_1+2\log h_2)(q)+24\pi \nabla_1 \beta(q,q)=0. $$
Thus if this function is never zero case three cannot happen either.  Theorem \ref{simple-2} is established. $\Box$

\medskip

\noindent{\bf Proof of Theorem \ref{simple-3}:}  First we observe that case two does not occur, since $\lambda_1^k$ is not greater than $2\lambda_2^k$. Next we need to rule out case three. The assumption $\lambda_1^k$ and $\lambda_1^k$ implies that $\mu_2^k>c\log \delta_k^{-1}$ for some $c>0$. Thus the $d_k$ terms can be removed from all related estimates.  The fact that $\lambda_2^k>\frac{3+\epsilon}4\lambda_1^k$ implies that $O(\delta_k^{-1}E)\to 0$ as $k\to \infty$. Thus at $q^k$ we have
$$\nabla (\log h_1+2\log h_2)(q^k)+24\pi\nabla_1\beta (q^k,q^k)=o(1), $$
a violation of the assumption that $\nabla (\log h_1+2\log h_2)+24\pi\nabla_1 \beta(\cdot, \cdot) $ is never zero.  Theorem \ref{simple-3} is established. $\Box$

\section{Appendix} In this appendix, we explain how to choose $p_1,p_2,p_3$ to make the following matrix invertible:
$$M=\left(\begin{array}{ccc}
\frac{\partial U}{\partial \Lambda}(p_1) & \frac{\partial U}{\partial \Lambda}(p_2) & \frac{\partial U}{\partial \Lambda}(p_3)\\
\frac{\partial U}{\partial \xi_1}(p_1) & \frac{\partial U}{\partial \xi_1}(p_2) & \frac{\partial U}{\partial \xi_1}(p_3)\\
\frac{\partial U}{\partial \xi_2}(p_1) & \frac{\partial U}{\partial \xi_2}(p_2) & \frac{\partial U}{\partial \xi_2}(p_3)
\end{array}
\right)$$
where $U$ is a global solution of
$$\Delta U+|y|^2 e^{U}=0,\quad \mbox{in}\quad \mathbb R^2, \quad \int_{\mathbb R^2}e^U<\infty, $$
with the expression $z=y_1+\sqrt{-1}y_2$
$$U(z)=\log \frac{\Lambda}{(1+\frac{\Lambda}{32}|z^2-\xi|^2)^2},\quad \xi\in \mathbb C. $$

The invertibility of $M$ is  crucial to approximate a bubbling solutions by a sequence of global solutions. Since the global solution has no symmetry, most of the traditional methods fail miserably. By adjusting parameters of global solutions infinitesimally we can make global solutions and bubbling solutions agree on three points (for $SU(3)$ Toda system) so that sharp estimates can be deduced from non-degeneracy results proved by Lin-Wei-Ye\cite{lwy}. In this appendix we employ the ideas in \cite{lwz-jems} to $SU(3)$, which may help readers understand better the more general situations discussed in that article.

Direct computation shows
\begin{align*}
\frac{\partial U}{\partial \Lambda}=\frac{1}{\Lambda}-\frac{1}{16}\frac{|z^{2}-\xi|^2}{1+\frac{\Lambda}{32}|z^{2}-\xi|^2},\\
\frac{\partial U}{\partial \xi}=\frac{\Lambda}{16}\frac{\bar z^{2}-\bar \xi}{1+\frac{\Lambda}{32}|z^{2}-\xi|^2},\nonumber \\
\frac{\partial U}{\partial \bar \xi}=\frac{\Lambda}{16}\frac{ z^{2}-\xi}{1+\frac{\Lambda}{32}|z^{2}-\xi|^2},\nonumber
\end{align*}

Obviously $M$ is invertible if and only if the following matrix is invertible:
$$M_1=\left(\begin{array}{ccc}
\frac{\partial U}{\partial \Lambda}(p_1) & \frac{\partial U}{\partial \Lambda}(p_2) & \frac{\partial U}{\partial \Lambda}(p_3)\\
\frac{\partial U}{\partial \xi}(p_1) & \frac{\partial U}{\partial \xi}(p_2) & \frac{\partial U}{\partial \xi}(p_3)\\
\frac{\partial U}{\partial \bar \xi}(p_1) & \frac{\partial U}{\partial \bar \xi}(p_2) & \frac{\partial U}{\partial \bar \xi}(p_3)
\end{array}
\right)$$
we choose $p_l$ to be
$$p_l=s^{1+\epsilon l}  e^{\sqrt{-1}\theta_l},\quad l=1,2,3 $$
where $s\ge \ge 1\ge \ge \epsilon>0 $ are constants to be determined later. Using crude expansion, we can write $M_1$ as
$$\left(\begin{array}{ccc}
\frac{1}{\Lambda}+O(s^{-2-2\epsilon}) & \frac{1}{\Lambda}+O(s^{-2-4\epsilon}) & \frac{1}{\Lambda}+O(s^{-2-6\epsilon}) \\
\frac{\Lambda}{16}s^{-2-2\epsilon}e^{-2i\theta_1}+O(s^{-4-4\epsilon}) & \frac{\Lambda}{16}s^{-2-4\epsilon}e^{-2i\theta_2}+O(s^{-4-8\epsilon}) &
\frac{\Lambda}{16}s^{-2-6\epsilon}e^{-2i\theta_3}+O(s^{-4-12\epsilon}) \\
\frac{\Lambda}{16}s^{-2-2\epsilon}e^{2i\theta_1}+O(s^{-4-4\epsilon}) & \frac{\Lambda}{16}s^{-2-4\epsilon}e^{2i\theta_2}+O(s^{-4-8\epsilon}) &
\frac{\Lambda}{16}s^{-2-6\epsilon}e^{2i\theta_3}+O(s^{-4-12\epsilon})
\end{array}
\right)
$$
Multiplying the first row by $\Lambda$, the second row and the third row by $\frac{16}{\Lambda}s^{2+6\epsilon}$, we change $M_1$ to
$$\left(\begin{array}{ccc}
1+O(s^{-2-2\epsilon}) & 1+O(s^{-2-4\epsilon}) & 1+O(s^{-2-6\epsilon}) \\
s^{4\epsilon}e^{-2i\theta_1}+O(s^{-2+2\epsilon}) & s^{2\epsilon}e^{-2i\theta_2}+O(s^{-2-2\epsilon}) &
e^{-2i\theta_3}+O(s^{-2-6\epsilon}) \\
s^{4\epsilon}e^{2i\theta_1}+O(s^{-2+2\epsilon}) & s^{2\epsilon}e^{2i\theta_2}+O(s^{-2-2\epsilon}) &
e^{2i\theta_3}+O(s^{-2-6\epsilon})
\end{array}
\right)
$$
By choosing $s>>1$ and $0<\epsilon<<1$ it is easy to see that the determinant of the matrix above is not zero if any only if the following matrix is non-zero:
$$M_2:=\left(\begin{array}{ccc}
1 & 1 & 1 \\
s^{4\epsilon}e^{-2i\theta_1} & s^{2\epsilon}e^{-2i\theta_2} &
e^{-2i\theta_3} \\
s^{4\epsilon}e^{2i\theta_1} & s^{2\epsilon}e^{2i\theta_2} &
e^{2i\theta_3}
\end{array}
\right)
$$
The determinant of $M_2$ is $2i \sin(2(\theta_2-\theta_1)) s^{6\epsilon}+O(s^{4\epsilon})$. Thus by choosing $s$ large, $\epsilon$ small and $\theta_1,\theta_2$ appropriately we can make $M_2$, as well as $M$, invertible.

\end{document}